\documentclass[reqno]{amsart}
\usepackage{hyperref}

\begin{document}
\title[On the constructive solution of an inverse Sturm-Liouville problem]{On the constructive solution of an inverse Sturm-Liouville problem}

\author{Avetik Pahlevanyan}

\address{Avetik Pahlevanyan \newline
Institute of Mathematics of National Academy of Sciences, 24/5 Baghramian ave., 0019, Yerevan, Armenia}
\email{apahlevanyan@instmath.sci.am}


\begin{abstract}
The necessary and sufficient conditions are found for the two sequences $\left\{\mu_n \right\}_{n=0}^{\infty}$ and $\left\{a_n \right\}_{n=0}^{\infty}$ to be the spectrum and the norming constants respectively, for a boundary value problem $-y'' + q \left( x \right) y = \mu y,$ $y \left(0\right)=0,$ $y\left( \pi \right)\cos \beta + y'\left( \pi \right)\sin \beta  = 0,$ $\beta \in \left( 0, \pi \right),$ with $q \in L_{\mathbb{R}}^1 \left[0, \pi \right].$
\end{abstract}

\maketitle
\numberwithin{equation}{section}
\newtheorem{theorem}{Theorem}[section]
\newtheorem{lemma}[theorem]{Lemma}
\newtheorem{definition}[theorem]{Definition}
\newtheorem{remark}[theorem]{Remark}

\section{Introduction and Statement of the Results}
\label{sec1}

Let us denote by $L\left(q, \alpha, \beta \right)$ the following Sturm-Liouville boundary value problem
\begin{equation}\label{eq1.1}
- y'' + q\left( x \right)y = \mu y \equiv \lambda^2 y, \; x \in \left( 0, \pi \right), \; \mu \in \mathbb{C},
\end{equation}
\begin{equation}\label{eq1.2}
y\left( 0 \right)\cos \alpha  + y'\left( 0 \right)\sin \alpha  = 0, \; \alpha \in \left( 0, \pi \right],
\end{equation}
\begin{equation}\label{eq1.3}
y\left( \pi \right)\cos \beta  + y'\left( \pi \right)\sin \beta  = 0, \; \beta \in \left[ 0, \pi \right),
\end{equation}
where $q$ is a real-valued, summable function on $\left[0, \pi \right]$ (we write $q \in L_{\mathbb{R}}^1\left[ 0, \pi \right]$). By $L\left(q,\alpha, \beta \right)$ we also denote the self-adjoint operator generated by the problem \eqref{eq1.1}--\eqref{eq1.3} in Hilbert space $L^2\left[0, \pi \right]$ (see \cite{Naimark:1969,Levitan_Sargsyan:1970}). It is well-known, that the spectrum of $L\left(q, \alpha, \beta \right)$ is discrete and consists of real, simple eigenvalues (see \cite{Naimark:1969,Levitan_Sargsyan:1970,Marchenko:1977}), which we denote by ${\mu}_n \left(q, \alpha, \beta \right),$ $n=0,1,2,\dots,$ emphasizing the dependence on $q,$ $\alpha$ and $\beta.$ We assume that eigenvalues $\mu_n$ are enumerated in increasing order:
\begin{equation*}
\mu_0 \left(q, \alpha, \beta \right) < \mu_1 \left(q, \alpha, \beta \right) < \dots < \mu_n \left(q, \alpha, \beta \right) < \dots \; .
\end{equation*}
Let $\varphi \left(x, \mu, \alpha \right)$ and $\psi \left(x, \mu, \beta \right)$ be the solutions of \eqref{eq1.1}, satisfying the initial conditions
\begin{equation}\label{eq1.4}
\varphi \left(0, \mu, \alpha \right)=\sin \alpha, \; \varphi' \left(0, \mu, \alpha \right)=-\cos \alpha,
\end{equation}
\begin{equation}\label{eq1.5}
\psi \left(\pi, \mu, \beta \right)=\sin \beta, \; \psi' \left(\pi, \mu, \beta \right)=-\cos \beta.
\end{equation}
It is well-known (\cite{Naimark:1969,Levitan_Sargsyan:1970,Harutyunyan_Hovsepyan:2005}) that for fixed $x,$ the functions $\varphi,$ $\varphi',$ $\psi,$ $\psi'$ are entire with respect to $\mu.$
We set by $W_{\alpha, \beta}\left(x, \mu \right)$ Wronskian of the solutions $\varphi \left(x, \mu, \alpha \right)$ and $\psi \left(x, \mu, \beta \right):$
\begin{equation}\label{eq1.6}
W_{\alpha, \beta}\left(x, \mu \right):=\varphi \left(x, \mu, \alpha \right) \psi' \left(x, \mu, \beta \right) - \varphi' \left(x, \mu, \alpha \right) \psi \left(x, \mu, \beta \right).
\end{equation}
By virtue of Liouville formula for the Wronskian (see, e.g., \cite{Coddington_Levinson:1955}) $W_{\alpha, \beta}\left(x, \mu \right)$ does not depend on $x,$ i.e.
\begin{equation}\label{eq1.7}
W_{\alpha, \beta}\left(\mu \right)=W_{\alpha, \beta}\left(x, \mu\right)=W_{\alpha, \beta}\left( \pi, \mu \right)=W_{\alpha, \beta}\left( 0, \mu \right).
\end{equation}
The eigenvalues $\mu_n={\mu_n}\left(q, \alpha, \beta \right),$ $n=0,1,2,\dots,$ of $L\left(q, \alpha, \beta \right)$ are the solutions of the equation $W_{\alpha, \beta}\left(\mu \right)=0$ (see, e.g., \cite{Marchenko:1977}). It is easy to see that the functions ${\varphi_n}\left(x\right):=\varphi(x, \mu_n, \alpha)$ and ${\psi_n}\left(x\right):=\psi\left(x, \mu_n, \beta\right),$ $n=0,1,2,\dots,$ are the eigenfunctions corresponding to the eigenvalue $\mu_n.$ Since the eigenvalues are simple, then there exist the numbers $\beta_n=\beta_n\left(q, \alpha, \beta\right),$ $n=0,1,2,\dots,$ such that
\begin{equation}\label{eq1.8}
\psi_n\left(x\right)=\beta_n \varphi_n\left(x\right), \beta_n \neq 0.
\end{equation}
The squares of the $L^2$ norms of these eigenfunctions:
\begin{equation}\label{eq1.9}
a_n=a_n \left(q, \alpha, \beta \right)=\int\limits_0^\pi {\varphi_n^2 \left(x\right)dx}, \quad b_n=b_n \left(q, \alpha, \beta \right)=\int\limits_0^\pi {\psi_n^2 \left(x\right)dx},
\end{equation}
are called norming constants.

In the works \cite{Chudov:1949,Marchenko:1952,Gelfand_Levitan:1951,Gasymov_Levitan:1964,Zhikov:1967,Isaacson_Trubowitz:1983,
Isaacson_Mckean_Trubowitz:1984,Poschel_Trubowitz:1987,Levitan:1984,Freiling_Yurko:2001} various aspects of the inverse Sturm-Liouville problem, consisting of recovering operator by two sequences, has been considered. Particularly, the following question has been studied:

"What kind must be the sequences $\left\{\mu_n\right\}_{n=0}^{\infty}$ and $\left\{a_n\right\}_{n=0}^{\infty},$ to be the spectrum and the norming constants of a problem $L \left(q, \alpha, \beta \right),$ respectively".

This question is well studied for the cases $q \in L_{\mathbb{R}}^2\left[ 0, \pi \right],$ $\alpha, \beta \in \left(0, \pi \right)$ and $\alpha=\pi,$ $\beta=0.$ For the case $q \in L_{\mathbb{R}}^1\left[ 0, \pi \right]$ and $\sin \alpha=0$ $\left(\alpha=\pi\right),$ $\beta \in \left(0, \pi\right)$ (analogously for $\alpha \in \left(0, \pi\right),$ $\sin \beta=0$ $\left(\beta=0\right)$) only some aspects of this question have been studied by the above-mentioned and the other authors. From the other point of view, an inverse problem for $L\left(q, \pi, \beta \right)$ (with $q \in L_{\mathbb{R}}^2\left[ 0, \pi \right]$ and $\beta \in \left(0, \pi\right)$) was considered in \cite{Dahlberg_Trubozitz:1984,Korotyaev_Chelkak:2009}.

To our knowledge, so far, the necessary and sufficient conditions for the sequences $\left\{\mu_n\right\}_{n=0}^{\infty}$ and $\left\{a_n\right\}_{n=0}^{\infty}$ to be the spectrum and the norming constants for the problem $L \left(q, \pi, \beta \right)$ with $q \in L_{\mathbb{R}}^1\left[0, \pi \right]$ (analogously for $L \left(q, \alpha, 0 \right)$) have not been found, which in particularly means that the constructive solution of the inverse Sturm-Liouville problem in this case is not given. Our goal is to find these conditions.

To formulate the main theorem of the work we give a notation.

In the paper \cite{Harutyunyan:2008} Harutyunyan introduced the concept of the function of $\delta_n \left(\alpha, \beta \right),$ which is defined as $\sqrt {\mu_n \left(0, \alpha, \beta \right)}-n:=\delta_n \left(\alpha, \beta \right),$ $n \geq 2$ and proved that $-1 \leq \delta_{n} \left(\alpha, \beta \right) \leq 1$ and it is a solution of the following transcendental equation:
\begin{multline}\label{eq1.10}
\delta_n \left(\alpha, \beta \right)=\dfrac{1}{\pi}\arccos \dfrac{\cos \alpha}{\sqrt {\left(n+{\delta_n} \left(\alpha, \beta \right)\right)^2 \sin^2 \alpha + \cos^2 \alpha}}- \\
-\dfrac{1}{\pi}\arccos \dfrac{\cos \beta}{\sqrt {\left(n+\delta_n \left(\alpha, \beta \right)\right)^2 \sin^2 \beta+\cos^2 \beta}}.
\end{multline}
Observe that, since $\arccos$ is a decreasing function, then the transcendental equation \eqref{eq1.10} has a unique solution for $\alpha=\pi$ and $\beta \in \left[0, \pi \right),$ that is why using \eqref{eq1.10} for the determination of $\delta_n \left(\pi, \beta \right)$ (see, e.g., \eqref{eq3.8}) is correct.

The main result of this paper is the following theorem:

\begin{theorem}\label{thm1.1}
For the two sequences $\left\{\mu_n\right\}_{n=0}^{\infty}$ and $\left\{a_n\right\}_{n=0}^{\infty}$ to be the spectrum and the norming constants of a problem $L \left(q, \pi, \beta \right),$ with $q \in L_{\mathbb{R}}^1 \left[0, \pi \right]$ and some $\beta \in \left(0, \pi\right),$ it is necessary and sufficient that the following relations hold:
\begin{equation}\label{eq1.11}
\sqrt {\mu_n} \equiv \lambda_n=n+\delta_n \left(\pi, \beta\right)+ \dfrac{c}{2\left(n+\delta_n\left(\pi, \beta\right)\right)}+l_n, \; \mu_n \neq \mu_m \, \left(n \neq m \right),
\end{equation}
\begin{equation}\label{eq1.12}
a_n =\dfrac{\pi}{2 {\left(n + \delta_n \left(\pi, \beta \right)\right)}^2} \left( 1+ \dfrac {2 \, s_n}{\pi \left(n + \delta_n \left( \pi, \beta \right)\right)}\right), \; a_n>0,
\end{equation}
where $c$ is a constant, the reminders $l_n=o\left(\dfrac{1}{n}\right)$ and $s_n=o\left(1\right)$ (when $n \to \infty$) are such that, the functions
\begin{equation}\label{eq1.13}
l\left(t\right)=\displaystyle\sum_{n=2}^{\infty} l_n \sin \left(n+\delta_n \left(\pi, \beta\right)\right)t
\end{equation}
and
\begin{equation}\label{eq1.14}
s\left(t\right) = \displaystyle \sum_{n=2}^{\infty} \dfrac {s_n}{n + \delta_n \left( \pi, \beta \right)} \cos \left(n + \delta_n \left( \pi, \beta \right)\right)t
\end{equation}
are absolutely continuous on arbitrary segment $\left[a, b \right] \subset \left(0, 2 \pi \right)$ (we will write
$l, s \in AC\left(0, 2\pi \right)$).
\end{theorem}

The outline of this paper is the following: in section \ref{sec2} we prove new expansion theorems, that are used in section \ref{sec3} for the proof of the asymptotics of the eigenvalues and the norming constants (this case needs more fine-grained assessment of the remainders) and in section \ref{sec4} for derivation of the analogue of Gelfand-Levitan equation for our case $\left(\alpha=\pi, \, q \in L_{\mathbb{R}}^1 \left[0,\pi \right] \right).$ In section \ref{sec5} the existence and uniqueness of the solution of this equation as well as reconstruction of the function $q$ (i.e. the reconstruction of differential equation \eqref{eq1.1}) and parameter $\tilde \beta$ (see Remark \ref{rem5.7}) are given. In section \ref{sec6} by taking two certain sequences and using reconstruction technique we find the potential $q$ and the parameter $\tilde \beta.$ Appendix (section \ref{app}) includes two auxiliary lemmas, which are used to obtain some results of the paper.

\section{Expansion theorems}
\label{sec2}

Completeness and expansion theorems for eigenfunctions of the Sturm-Liouville boundary-value problem have been proven since XIX-th century. One of the main theorems of the spectral theory of differential operators is as follows (see \cite{Naimark:1969}):

\begin{theorem}(\cite[p. 90]{Naimark:1969})\label{thm2.1}
Every function in the domain of self-adjoint differential operator can be expanded as a uniformly convergent generalized Fourier series for eigenfunctions of this operator.
\end{theorem}

This result cannot be applied for the functions which do not belong to the domain of self-adjoint differential operator. On the other hand, it was proved that absolutely continuity of the function $f$ on $\left[0,\pi\right]$ is sufficient for Fourier series for the eigenfunctions of the Sturm-Liouville operator $L \left( q, \alpha, \beta \right),$ $q \in L_{\mathbb{R}}^2 \left[0, \pi \right],$ $\alpha, \beta \in \left(0, \pi \right)$ to converge uniformly to $f$ (see \cite{Coddington_Levinson:1955,Levitan_Sargsyan:1970,Levitan_Sargsyan:1988,Freiling_Yurko:2001}):

\begin{theorem}(\cite{Freiling_Yurko:2001})\label{thm2.2}
Let $q \in L_{\mathbb{R}}^2 \left[0, \pi \right],$ $\alpha, \beta \in \left(0, \pi \right)$ and $f$ be an absolutely continuous function on $\left[0, \pi \right].$ Then
$$ \mathop {\lim}\limits_{N \to \infty}\mathop {\operatorname{max}}\limits_{x \in \left[0, \pi\right]}\left|f\left(x\right)-\sum\limits_{n = 0}^{N} c_n {\varphi_n \left(x \right)}\right|=0, \; c_n=\dfrac{1}{a_n}{\displaystyle\int_{0}^{\pi}f\left(t\right){\varphi_n \left(t \right)}dt},$$
where $\varphi_n \left( x \right) \equiv \varphi \left(x, \mu_n \left(q, \alpha, \beta \right), \alpha \right).$
\end{theorem}

We prove that the analogous results for the problems $L \left( q, \pi, \beta \right),$ $\beta \in \left(0,\pi \right)$ and $L \left( q, \alpha, 0 \right),$ $\alpha \in \left(0,\pi\right)$ are also true:

\begin{theorem}\label{thm2.3}
Let $q \in L_{\mathbb{R}}^1\left[ {0,\pi } \right],$ $\alpha=\pi,$ $\beta \in \left(0, \pi \right)$ and $f$ be an absolutely continuous function on $\left[0, \pi\right].$ Then for arbitrary $a \in \left(0, \pi \right)$
\begin{equation}\label{eq2.1}
\mathop {\lim}\limits_{N \to \infty}\mathop {\operatorname{max}}\limits_{x \in \left[a, \pi\right]}\left|f\left(x\right)-\sum\limits_{n = 0}^{N} c_n {\varphi_n \left(x \right)}\right|=0, \; c_n=\dfrac{1}{a_n}{\displaystyle\int_{0}^{\pi}f\left(t\right){\varphi_n \left(t \right)}dt},
\end{equation}
where $\varphi_n \left( x \right) \equiv \varphi \left(x, \mu_n \left(q, \pi, \beta \right), \pi \right) \equiv \varphi \left(x, \mu_n, \pi \right).$
\end{theorem}

\begin{theorem}\label{thm2.4}
Let $q \in L_{\mathbb{R}}^1\left[ {0,\pi } \right],$ $\alpha \in \left(0, \pi \right),$ $\beta=0$ and  $f$ be an absolutely continuous function on $\left[0, \pi\right].$ Then for arbitrary $b \in \left(0, \pi \right)$
\begin{equation}\label{eq2.2}
\mathop {\lim}\limits_{N \to \infty}\mathop {\operatorname{max}}\limits_{x \in \left[0, b \right]}\left|f\left(x\right)-\sum\limits_{n = 0}^{N} c_n {\varphi_n \left(x \right)}\right|=0, \; c_n=\dfrac{1}{a_n}{\displaystyle\int_{0}^{\pi}f\left(t\right){\varphi_n \left(t \right)}dt},
\end{equation}
where $\varphi_n \left( x \right) \equiv \varphi \left(x, \mu_n \left(q, \alpha, 0 \right), \alpha \right).$
\end{theorem}

We provide the proof for Theorem \ref{thm2.3}. Theorem \ref{thm2.4} can be proved similarly. We are going to adapt the proof of the Theorem \ref{thm2.2} by Cauchy's contour integral method.

\begin{proof}
For $\left|\lambda\right| \to \infty,$ the following asymptotic formulae hold (\cite{Naimark:1969,Harutyunyan_Hovsepyan:2005,Marchenko:1952,Chudov:1949,Atkinson:1964,Harutyunyan:2016})
\begin{equation}\label{eq2.3}
\varphi \left(x, \mu, \pi \right) := \varphi_{\pi} \left(x, \mu \right) \equiv \varphi_{\pi} \left(x, \lambda^2 \right)=\dfrac{\sin \lambda x}{\lambda}+O\left(\dfrac{e^{\left|Im \lambda \right|x}}{\left|\lambda \right|^2}\right),
\end{equation}
\begin{equation}\label{eq2.4}
\varphi' \left(x, \mu, \pi \right) := \varphi'_{\pi} \left(x, \mu \right) \equiv \varphi'_{\pi} \left(x, \lambda^2 \right)=\cos \lambda x+O\left(\dfrac{e^{\left|Im \lambda \right|x}}{\left|\lambda \right|}\right),
\end{equation}
\begin{multline}\label{eq2.5}
\psi \left(x, \mu, \beta \right) := \psi_{\beta} \left(x, \mu \right) \equiv \psi_{\beta} \left(x, \lambda^2 \right)=\cos \lambda \left( \pi-x \right) \sin \beta + \dfrac{\sin \lambda \left( \pi-x \right)}{\lambda}\cos \beta + \\
+O\left(\dfrac{e^{\left|Im \lambda \right|(\pi-x)}}{\left|\lambda \right|}\right)\sin \beta + O\left(\dfrac{e^{\left|Im \lambda \right|(\pi-x)}}{\left|\lambda \right|^2}\right)\cos \beta,
\end{multline}
\begin{multline}\label{eq2.6}
\psi' \left(x, \mu, \beta \right) := \psi'_{\beta} \left(x, \mu \right) \equiv \psi'_{\beta} \left(x, \lambda^2 \right)=\left( \lambda \sin \lambda \left( \pi-x \right) + O\left(e^{\left|Im \lambda \right|(\pi-x)} \right) \right)\sin \beta - \\
- \left( \cos \lambda \left( \pi-x \right) + O\left(\dfrac{e^{\left|Im \lambda \right|(\pi-x)}}{\left|\lambda \right|}\right)\right) \cos \beta.
\end{multline}
From \eqref{eq1.6}, \eqref{eq1.7} and \eqref{eq2.5} for Wronskian $W_{\pi, \beta}\left(\mu \right)$ we have the following estimates
\begin{multline}\label{eq2.7}
W_{\pi, \beta}\left(\mu \right) \equiv W_{\pi, \beta}\left(\lambda^2 \right) = -\psi_{\beta} \left(0, \mu \right) = -\cos \lambda \pi \sin \beta - \dfrac{\sin \lambda \pi}{\lambda}\cos \beta + \\
+O\left(\dfrac{e^{\left|Im \lambda \right|\pi}}{\left|\lambda \right|}\right)\sin \beta + O\left(\dfrac{e^{\left|Im \lambda \right|\pi}}{\left|\lambda \right|^2}\right)\cos \beta.
\end{multline}
Denote by ${\mathbb{Z}}_{1/6}$ the following domain of the complex plane $\mathbb{C}$:
$${\mathbb{Z}}_{1/6}=\left\{\lambda \in \mathbb{C}: \; \left|\lambda-\dfrac{n}{2}\right| \geq \dfrac{1}{6}, \; n \in \mathbb{Z} \right\}.$$
The following lemma have been proven in \cite{Harutyunyan:2010} by the methods, which had used in \cite[Lemma 1 on p. 27]{Poschel_Trubowitz:1987} (see, also \cite{Marchenko:1977,Freiling_Yurko:2001})
\begin{lemma}(\cite{Harutyunyan:2010})\label{lem2.5}
If $\lambda \in {\mathbb{Z}}_{1/6}$, then
\begin{equation}\label{eq2.8}
\left|\sin \pi \lambda \right| \geq \dfrac{1}{7}e^{\left|Im\lambda\right|\pi}, \; \left|\cos \pi \lambda \right| \geq \dfrac{1}{7}e^{\left|Im\lambda\right|\pi}.
\end{equation}
\end{lemma}
It follows from \eqref{eq2.7} and \eqref{eq2.8} that for sufficiently large $\lambda^{*}>0,$  there is a constant $C_1>0$ such that
\begin{equation}\label{eq2.9}
\left|W_{\pi, \beta} \left( \lambda^2 \right)\right| \geq C_1 e^{\left|Im \lambda \right|\pi}\sin \beta, \; \mbox{when} \; \lambda \in {\mathbb{Z}}_{1/6}, \; \left|\lambda \right| > \lambda^{*}.
\end{equation}
Let us consider the following boundary value problem
\begin{equation}\label{eq2.10}
- y'' + q\left( x \right)y = \mu y-f \left( x \right), \; x \in \left( {0,\pi } \right), \; \mu \in \mathbb{C}, \; f \in L^1 \left[0,\pi\right],
\end{equation}
\begin{equation}\label{eq2.11}
y\left( 0 \right)=0, \; y\left( \pi \right)\cos \beta  + y'\left( \pi \right)\sin \beta  = 0, \; \beta  \in \left( {0,\pi } \right).
\end{equation}
It is well-known and can be easily verified that the solution $y\left(x, \mu, f \right)$ of the boundary value problem \eqref{eq2.10}--\eqref{eq2.11} can be written in the following form (see, for example, \cite{Naimark:1969,Levitan_Sargsyan:1988})
\begin{multline}\label{eq2.12}
y \left(x, \mu, f \right) = \dfrac{1}{W_{\pi, \beta}\left( \mu \right)} \psi_{\beta}\left( x, \mu \right) \displaystyle\int_{0}^{x}f\left(t\right)\varphi_{\pi} \left(t, \mu \right)dt+\\
+\dfrac{1}{W_{\pi, \beta}\left( \mu \right)}\varphi_{\pi} \left(x, \mu \right) \displaystyle\int_{x}^{\pi}f\left(t\right)\psi_{\beta} \left(t, \mu \right)dt.
\end{multline}
Since $\varphi,$ $\psi$ and $W_{\pi, \beta}$ are entire functions of $\mu$, then we see that $y \left(x, \mu, f \right)$ is a meromorphic function of $\mu,$ with poles in the zeros of $W_{\pi, \beta}$ or, that is the same, in eigenvalues $\mu_{n}, n=0,1,2,\dots.$ Since $\dot{W}_{\pi, \beta}\left(\mu_n\right) \equiv \dfrac{d}{d\mu}W_{\pi, \beta} \left(\mu_n\right)=\beta_n a_n$ (see \cite[Lemma 1.1.1]{Freiling_Yurko:2001}), then using \eqref{eq1.8}, we get the residue
\begin{multline}\label{eq2.13}
\mathop {\operatorname{Res}}\limits_{\mu  = {\mu _n}} {y \left(x, \mu, f \right)} = \dfrac{1}{\dot{W}_{\pi, \beta}\left( \mu_n \right)} \psi_{\beta}\left( x, \mu_n \right) \displaystyle\int_{0}^{x}f\left(t\right)\varphi_{\pi} \left(t, \mu_n \right)dt+\\
+\dfrac{1}{\dot{W}_{\pi, \beta}\left( \mu_n \right)}\varphi_{\pi} \left(x, \mu_n \right) \displaystyle\int_{x}^{\pi}f\left(t\right)\psi_{\beta} \left(t, \mu_n \right)dt=\dfrac{\beta_n}{\dot{W}_{\pi, \beta}\left( \mu_n \right)}\varphi_{\pi}\left( x, \mu_n \right) \displaystyle\int_{0}^{\pi}f\left(t\right)\varphi_{\pi} \left(t, \mu_n \right)dt=\\
=\dfrac{1}{a_n}\varphi_\pi \left( x, \mu_n \right)\displaystyle\int_{0}^{\pi}f\left( t \right) \varphi_\pi \left(t, \mu_n \right)dt.
\end{multline}
It follows from \eqref{eq2.3}, \eqref{eq2.5}, \eqref{eq2.9} and \eqref{eq2.12} that there are positive numbers $C,$ $C_2,$ $C_3,$ $C_4$ such that the following estimates hold for $\lambda \in {\mathbb{Z}}_{1/6},$ $\left|\lambda \right| > \lambda^{*}:$
\begin{multline}\label{eq2.14}
\left| y \left(x, \lambda^2, f \right) \right| \leq \dfrac{\left| \psi_{\beta} \left(x, \lambda^2 \right) \right|\mathop {\operatorname{max}}\limits_{t \in \left[0, x \right]}\left|\varphi_{\pi} \left(t, \lambda^2 \right) \right|\displaystyle\int_{0}^{x}\left|f\left( t \right)\right|dt}{C_1 e^{\left| Im \lambda \right|\pi}\sin \beta}+\\
+\dfrac{\left|\varphi_{\pi} \left(x, \lambda^2 \right) \right|\mathop {\operatorname{max}}\limits_{t \in \left[x, \pi \right]}\left|\psi_{\beta} \left(t, \lambda^2 \right) \right| \displaystyle\int_{x}^{\pi}\left|f\left( t \right)\right|dt}{C_1 e^{\left| Im \lambda \right|\pi} \sin \beta} \leq \\
\leq \dfrac{e^{\left| Im \lambda \right|(\pi-x)} \left(\sin \beta+\dfrac{ \left| \cos \beta \right| }{\left| \lambda \right|}+C_3\dfrac{ \sin \beta }{\left| \lambda \right|}+C_4\dfrac{ \left| \cos \beta \right| }{\left| \lambda \right|^2} \right)}{C_1 e^{\left| Im \lambda \right|\pi} \sin \beta} \times \\
\times e^{\left| Im \lambda \right|x} \left(\dfrac{1}{\left| \lambda \right|}+C_2\dfrac{1}{\left| \lambda \right|^2} \right)\displaystyle\int_{0}^{\pi}\left|f\left( t \right)\right|dt \leq \\
\leq \dfrac{1}{C_1} \displaystyle\int_{0}^{\pi}\left|f\left( t \right)\right|dt \left( \dfrac{1}{\left| \lambda \right|}+O\left( \dfrac{1}{\left| \lambda \right|^2} \right) \right) \leq \dfrac{C}{\left| \lambda \right|}.
\end{multline}
Let us now consider a function $f \in AC\left[0, \pi \right].$ Using the fact that $\varphi_{\pi} \left(x, \mu \right)$ and $\psi_{\beta} \left(x, \mu \right)$ are the solutions of \eqref{eq1.1}, we can rewrite the representation \eqref{eq2.12} for $y \left(x, \mu, f \right)$ in the following form (compare with \cite{Freiling_Yurko:2001}):
\begin{equation}\label{eq2.15}
y \left(x, \mu, f \right) = \dfrac{f\left(x\right)}{\mu}+f\left(0\right) \dfrac{\psi_{\beta}\left(x, \mu \right)}{\mu W_{\pi, \beta}\left(\mu\right)}+\dfrac{Z_1\left(x, \mu, \pi, \beta, f' \right)} {\mu}+\dfrac{Z_2\left(x, \mu, \pi, \beta \right)}{\mu},
\end{equation}
where
\begin{equation}\label{eq2.16}
Z_1\left(x, \mu, \pi, \beta, f' \right)=\dfrac{\psi_{\beta}\left( x, \mu \right) \displaystyle\int_{0}^{x} f'\left(t\right) \varphi'_{\pi} \left(t, \mu \right)dt +\varphi_{\pi}\left( x, \mu \right) \displaystyle\int_{x}^{\pi}f'\left(t\right)\psi'_{\beta} \left(t, \mu \right)dt}{W_{\pi, \beta}\left( \mu \right)},
\end{equation}
\begin{multline}\label{eq2.17}
Z_2\left(x, \mu, \pi, \beta \right)=-f\left(\pi\right)\psi'_{\beta} \left(\pi, \mu \right) \dfrac{\varphi_{\pi}\left( x, \mu \right)}{W_{\pi, \beta}\left( \mu \right)}+y \left(x, \mu, qf \right)= \\
=f\left(\pi\right)\cos \beta \dfrac{\varphi_{\pi}\left( x, \mu \right)}{W_{\pi, \beta}\left( \mu \right)}+y \left(x, \mu, qf \right).
\end{multline}
Let us show that
\begin{equation}\label{eq2.18}
\mathop {\lim }\limits_{\left| \lambda  \right| \to \infty \hfill \atop \lambda  \in {\mathbb{Z}}_{1/6}} \hfill \mathop {\operatorname{max}}\limits_{x \in \left[0, \pi\right]}\left|Z_1\left(x, \mu, \pi, \beta, f' \right)\right|=0.
\end{equation}
First, we suppose that $f'$ is absolutely continuous function on $\left[0, \pi\right].$ Then there exists $f'' \in L^1\left[ {0,\pi } \right]$ and \eqref{eq2.16} can be written in the following form
\begin{multline*}
Z_1\left(x, \mu, \pi, \beta, f' \right)=\dfrac{\psi_{\beta}\left( x, \mu \right)}{W_{\pi, \beta}\left( \mu \right)} \left(\left.\varphi_{\pi}\left(t, \mu\right)f'\left(t\right)\right|_{t=0}^{x}- \displaystyle\int_{0}^{x} f''\left(t\right) \varphi_{\pi} \left(t, \mu \right)dt\right)+\\
+\dfrac{\varphi_{\pi}\left(x, \mu \right)}{W_{\pi, \beta}\left( \mu \right)} \left(\left.\psi_{\beta}\left(t, \mu\right)f'\left(t\right)\right|_{t=x}^{\pi}- \displaystyle\int_{x}^{\pi} f''\left(t\right) \psi_{\beta} \left(t, \mu \right)dt\right)=\dfrac{\varphi_{\pi}\left(x, \mu \right)}{W_{\pi, \beta}\left( \mu \right)} f'\left(\pi \right) \sin \beta-\\
-\dfrac{\psi_{\beta}\left( x, \mu \right) \displaystyle\int_{0}^{x} f''\left(t\right) \varphi_{\pi} \left(t, \mu \right)dt +\varphi_{\pi}\left( x, \mu \right) \displaystyle\int_{x}^{\pi}f''\left(t\right)\psi_{\beta} \left(t, \mu \right)dt}{W_{\pi, \beta}\left( \mu \right)}.
\end{multline*}
By virtue of \eqref{eq2.3}--\eqref{eq2.6} and \eqref{eq2.9} we obtain that there is a number $C>0,$ such that
\begin{equation*}
\mathop {\operatorname{max}}\limits_{x \in \left[0, \pi\right]}\left|Z_1\left(x, \mu, \pi, \beta, f' \right)\right| \leq \dfrac{C}{\left| \lambda \right|}, \; \mbox{when} \; \lambda \in {\mathbb{Z}}_{1/6}, \; \left|\lambda \right| > \lambda^{*}.
\end{equation*}
It implies \eqref{eq2.18} in the case $f' \in AC\left[0,\pi\right].$

Now, let us turn to the general case $g := f' \in L^1 \left[0, \pi \right].$ Fix $\epsilon>0$ and choose an absolutely continuous  function $g_\epsilon,$ such that
\begin{equation*}
\displaystyle\int_{0}^{\pi} \left|g\left(t\right)-{g_\epsilon}\left(t\right)\right|dt<\dfrac{C_1 \sin \beta }{16} \, \epsilon.
\end{equation*}
Then, according to \eqref{eq2.3}--\eqref{eq2.6}, \eqref{eq2.9} and \eqref{eq2.16} for $\lambda \in {\mathbb{Z}}_{1/6}, \; \left|\lambda \right| > \lambda^{*},$ we have
\begin{multline*}
\mathop {\operatorname{max}}\limits_{x \in \left[0, \pi\right]}\left|Z_1\left(x, \mu, \pi, \beta, g \right)\right| \leq \mathop {\operatorname{max}}\limits_{x \in \left[0, \pi\right]}\left|Z_1\left(x, \mu, \pi, \beta, {g_\epsilon}\right)\right|+\mathop {\operatorname{max}}\limits_{x \in \left[0, \pi\right]}\left|Z_1\left(x, \mu, \pi, \beta, g-{g_\epsilon}\right)\right| \leq \\
\leq \dfrac{C(\epsilon)}{\left| \lambda \right|}+\dfrac{C_1 \sin \beta }{16} \, \epsilon \mathop {\operatorname{max}}\limits_{x \in \left[0, \pi\right]}\left(\dfrac{\left|\psi_{\beta}\left(x, \mu \right)\right|\mathop {\operatorname{max}}\limits_{t \in \left[0, x \right]}\left|\varphi'_{\pi}\left(t, \mu\right)\right|+\left|\varphi_{\pi}\left(x, \mu \right)\right|\mathop {\operatorname{max}}\limits_{t \in \left[0, x \right]}\left|\psi'_{\beta}\left(t, \mu\right)\right|}{C_1 e^{\left|Im \lambda \right|\pi}\sin \beta}\right)\leq \\
\leq \dfrac{C(\epsilon)}{\left| \lambda \right|}+\dfrac{C_1 \sin \beta }{16} \, \epsilon \mathop {\operatorname{max}}\limits_{x \in \left[0, \pi\right]}\left(\dfrac{8 e^{\left|Im \lambda \right|\pi}}{C_1 e^{\left|Im \lambda \right|\pi}\sin \beta}\right) \leq \dfrac{C(\epsilon)}{\left| \lambda \right|}+\dfrac{\epsilon}{2}.
\end{multline*}
It is easy to see, that if we choose $\lambda_\epsilon^{*}=\dfrac{2C\left(\epsilon\right)}{\epsilon},$ then for $\lambda \in {\mathbb{Z}}_{1/6},$ $\left|\lambda \right| > \lambda_\epsilon^{*}$ we have $\mathop {\operatorname{max}}\limits_{x \in \left[0, \pi\right]}\left|Z_1\left(x, \mu, \pi, \beta\right)\right| \leq \epsilon.$  Due to the arbitrariness $\epsilon>0,$ we arrive at \eqref{eq2.18}.

Now, we estimate $Z_2\left(x, \mu, \pi, \beta\right)$ (see \eqref{eq2.17}). Since $qf \in L^1\left[ {0,\pi } \right],$ then the estimates in \eqref{eq2.14} are also true for $y\left( x, \mu, qf \right).$ Using \eqref{eq2.3}, \eqref{eq2.9}, \eqref{eq2.14} and the fact that $\sin \beta \neq 0$ we get the following estimates for $\lambda \in {\mathbb{Z}}_{1/6}, \; \left|\lambda \right| > \lambda^{*}$
\begin{multline}\label{eq2.19}
\mathop {\operatorname{max}}\limits_{x \in \left[0, \pi\right]}\left|Z_2\left(x, \mu, \pi, \beta\right)\right| \leq \mathop {\operatorname{max}}\limits_{x \in \left[0, \pi\right]} \left|f\left(\pi\right)\cos \beta \dfrac{\varphi_{\pi}\left( x, \mu \right)}{W_{\pi, \beta}\left( \mu \right)}\right| + \mathop {\operatorname{max}}\limits_{x \in \left[0, \pi\right]}\left|y \left(x, \mu, qf \right)\right| \leq \\
\leq \left|f\left(\pi\right)\cos \beta \dfrac{C_5 e^{\left|Im \lambda\right|\pi}}{\left|\lambda\right|C_1 e^{\left|Im \lambda\right|\pi}\sin \beta}\right| +\dfrac{C_6}{\left|\lambda\right|} \leq \dfrac{C_5}{C_1}\dfrac{\left|f\left(\pi\right)\cot \beta\right|} {\left|\lambda\right|} +\dfrac{C_6}{\left|\lambda\right|} \leq \dfrac{C_7}{\left|\lambda\right|},
\end{multline}
where $C_5,$ $C_6,$ $C_7$ are positive numbers.

Consider the following contour integral
\begin{equation}\label{eq2.20}
I_{N}\left(x\right)=\dfrac{1}{2\pi i}\oint_{\Gamma_{N}}y \left(x, \mu, f \right)d \mu,
\end{equation}
where $\Gamma_{N}=\left\{\mu : \left|\mu\right|=\left(N+\dfrac{3}{4}\right)^2\right\}$ (with counterclockwise circuit).
On one hand, using the Cauchy's residue theorem (see \cite{Shabat:1985}), from \eqref{eq2.13} we get
\begin{equation}\label{eq2.21}
I_{N}\left(x\right)= \sum\limits_{n = 0}^{N} \dfrac{1}{a_n}\displaystyle\int_{0}^{\pi}f\left( t \right){\varphi_{\pi} \left(t, \mu_n \right)}dt \varphi_{\pi} \left(x, \mu_n \right).
\end{equation}
On the other hand, from \eqref{eq2.15}, \eqref{eq2.18} and \eqref{eq2.19} we obtain
\begin{equation}\label{eq2.22}
I_{N}\left(x\right)= f\left(x\right)+f(0)\dfrac{1}{2\pi i}\oint_{\Gamma_{N}} \dfrac{\psi_{\beta}\left(x, \mu \right)}{\mu W_{\pi, \beta}\left(\mu\right)}d \mu+{\epsilon_{N}}\left(x\right),
\end{equation}
where ${\epsilon_{N}}\left(x\right),$ according to \eqref{eq2.18} and \eqref{eq2.19}, uniformly converges to $0:$
$$ \mathop {\lim}\limits_{N \to \infty}\mathop {\operatorname{max}}\limits_{x \in \left[0, \pi\right]}\left|{\epsilon_{N}}\left(x\right)\right|=0.$$
Without loss of generality we assume that $\mu=0$ is not an eigenvalue of the problem $L\left(q,\pi,\beta\right).$ Indeed, if $0$ is an eigenvalue, then there is a number $c$ such that $\mu_n +c \neq 0,$ $n=0,1,2,\dots$ are the eigenvalues of the problem $L\left(q+c, \pi, \beta \right)$ with the same eigenfunctions $\varphi_n$ and norming constants $a_n$ as for $L\left(q, \pi, \beta \right).$ Then the function $\dfrac{\psi_{\beta}\left(x, \mu \right)}{\mu W_{\pi, \beta}\left(\mu\right)}$ has only first order poles and using Cauchy's residue theorem (see \cite{Shabat:1985}), we can easily calculate that
\begin{multline}\label{eq2.23}
\phi_N \left(x \right) :=\dfrac{1}{2\pi i}\oint_{\Gamma_{N}} \dfrac{\psi_{\beta}\left(x, \mu \right)}{\mu W_{\pi, \beta}\left(\mu\right)}d \mu=\mathop {\operatorname{Res}}\limits_{\mu=0}\dfrac{\psi_{\beta}\left(x, \mu \right)}{\mu W_{\pi, \beta}\left(\mu\right)} +\sum\limits_{n = 0}^{N}\mathop {\operatorname{Res}}\limits_{\mu=\mu_n}\dfrac{\psi_{\beta}\left(x, \mu \right)}{\mu W_{\pi, \beta}\left(\mu\right)}=\\
=\dfrac{\psi_{\beta}\left(x, 0 \right)}{W_{\pi, \beta}\left(0\right)}+\sum\limits_{n = 0}^{N} \dfrac{\psi_{\beta}\left(x, \mu_n \right)}{\mu_n \dot{W}_{\pi, \beta}\left(\mu_n\right)}= \dfrac{\psi_{\beta}\left(x, 0 \right)}{W_{\pi, \beta}\left(0\right)}+\sum\limits_{n = 0}^{N} \dfrac{\beta_n \varphi_{\pi}\left(x, \mu_n \right)}{\mu_n \beta_n a_n}=\\
=\dfrac{\psi_{\beta}\left(x, 0 \right)}{W_{\pi, \beta}\left(0\right)}+\sum\limits_{n = 0}^{N} \dfrac{1}{\mu_n a_n} {\varphi_n \left(x \right)}.
\end{multline}
Now let us show that the sequence $\phi_N \left(x \right)$ converges to $0$ (when $N \to \infty$) uniformly on segment $\left[a, \pi \right],$ for arbitrary $a \in \left(0, \pi \right).$

Since $\varphi_n\left(x\right)=\dfrac{\sin \left(n+\frac{1}{2}\right)x}{n+\frac{1}{2}}+O\left(\dfrac{1}{n^2}\right)$ uniformly on $\left[0, \pi \right]$  (see \eqref{eq2.3}), $\mu_n=\mu_n \left(q,\pi,\beta\right)=\left(n+\dfrac{1}{2}\right)^2+O\left(1\right)$ (see \cite[Theorem 1 on p. 286 and asymptotical estimate for $\delta_n \left(\pi, \beta \right)$ on p. 292]{Harutyunyan:2008}) and $a_n=a_n\left(q,\pi,\beta\right)=\dfrac{\pi}{2\left(n+\frac{1}{2}\right)^2} \left(1+o\left(\dfrac{1}{n}\right)\right)$ (see \cite[Theorem 1.1 on pp. 9--10]{Harutyunyan_Pahlevanyan:2016}), then $\phi_N \left(x \right)$ (see \eqref{eq2.23}) can be written in the following form:
\begin{equation*}
\phi_N \left(x \right)=\dfrac{\psi_{\beta}\left(x, 0 \right)}{W_{\pi, \beta}\left(0\right)}+\dfrac{2}{\pi}\sum\limits_{n = 0}^N \dfrac{\sin \left(n+\frac{1}{2}\right)x}{n+\frac{1}{2}}+\sum\limits_{n = 0}^N q_n \left(x\right),
\end{equation*}
where $q_n \left(x \right)=O\left(\dfrac{1}{n^2}\right)$ uniformly on $\left[0, \pi \right].$

Since $\displaystyle\sum\limits_{n = 0}^\infty \dfrac{\sin \left(n+\frac{1}{2}\right)x}{n+\frac{1}{2}}=\dfrac{\pi}{2},$ $0 < x < 2\pi$ (see, for example, \cite[formula 37 on p. 578]{Bronshtein_Semendyayev:1998}), then the sequence $\phi_N \left(x \right)$ converges to the continuous function $\phi \left(x \right)$ (when $N \to \infty$) uniformly on segment $\left[a, \pi \right],$ for arbitrary $a \in \left(0, \pi \right).$

Now, to prove that $\phi \left(x \right) \equiv 0, \, x \in \left(0, \pi \right],$ it is sufficient to prove that $\phi=0$ in $L^2 \left(0,\pi\right).$

To this aim, by doing some calculations we get:
\begin{equation}\label{eq2.24}
\int\limits_0^\pi \phi \left(x \right) \varphi_m \left(x \right)dx = \dfrac{1}{W_{\pi, \beta}\left(0\right)} \int\limits_0^\pi \psi_{\beta}\left(x, 0 \right) \varphi_m \left(x\right)dx+\dfrac{1}{\mu_m}, \, m=0,1,2,\dots
\end{equation}
and
\begin{multline}\label{eq2.25}
\mu_m \int\limits_0^\pi \psi_{\beta}\left(x, 0 \right) \varphi_m \left(x\right)dx=\int\limits_0^\pi \left(\varphi_m \left(x\right)\psi''_{\beta}\left(x, 0 \right)-\varphi''_m \left(x\right) \psi_{\beta}\left(x, 0 \right)\right)dx=\\
=\left.\left(\varphi_m \left(x\right)\psi'_{\beta}\left(x, 0 \right)-\varphi'_m \left(x\right) \psi_{\beta}\left(x, 0 \right)\right)\right|_{0}^{\pi}=\psi_{\beta}\left(0, 0 \right)=-W_{\pi, \beta}\left(0\right).
\end{multline}

It follows from \eqref{eq2.24} and \eqref{eq2.25} that
\begin{equation*}
\int\limits_0^\pi \phi \left(x \right) \varphi_m \left(x \right)dx =0, \, m=0,1,2,\dots.
\end{equation*}

Since the system of eigenfunctions $\left\{\varphi_{m} \left(x \right)\right\}_{m=0}^\infty$ of the boundary value problem $L(q, \pi, \beta)$ is complete in $L^2 \left(0, \pi \right),$ then $\phi=0$ in $L^2 \left(0,\pi\right).$

Comparing this result with \eqref{eq2.21}, \eqref{eq2.22} and passing to the limit $\left(N \to \infty \right)$ we arrive at \eqref{eq2.1}. Theorem \ref{thm2.3} is proved.
\end{proof}

The attentive reader might offer the equiconvergence argument to obtain the result. Below in remark \ref{rem2.6} we briefly explain why this argument can't be applied in our case.

\begin{remark}\label{rem2.6}
It is well known that one of the proofs of the classical Theorem \ref{thm2.2} is based on the so-called uniformly equiconvergence theorem, which states that the expansion in eigenfunctions of the problem $L \left(q, \alpha, \beta \right),$ $\alpha, \beta \in \left(0,\pi \right)$ is equivalent to the expansion in eigenfunctions of the problem $L \left(0, \dfrac{\pi}{2}, \dfrac{\pi}{2} \right),$ i.e., $\left\{\cos nx\right\}_{n \geq 0}$ (see \cite{Coddington_Levinson:1955,Levitan_Sargsyan:1970,Levitan_Sargsyan:1988}). Further, the Dirichlet-Jordan theorem (see \cite[pp. 121--122]{Bari:1961}) can be applied and Theorem \ref{thm2.2} will be proved. On one hand, it is easy to set up a similar equiconvergence assertion for expansion in eigenfunctions of the problem $L \left(q, \pi, \beta \right),$ $\beta \in \left(0, \pi \right).$ This is equivalent to the expansion in $\sin \left(n+\dfrac{1}{2}\right)x,$ $n=0,1,2,\dots,$ i.e., eigenfunctions of the problem $L \left(0, \pi, \dfrac{\pi}{2} \right)$ (see \cite[remark on p. 304]{Coddington_Levinson:1955} and \cite[remark on p. 71]{Levitan_Sargsyan:1970}). On the other hand, to the best of our knowledge, there are no analogous result of the Dirichlet-Jordan theorem for expansion in the system of functions $\left\{\sin \left(n+\dfrac{1}{2}\right)x\right\}_{n \geq 0}$ (for the nearest result see \cite[Theorem 2.6]{Iserles_Norsett:2008}). To overcome this difficulty, we have proved Theorem \ref{thm2.3} using the Cauchy's contour integral method.
\end{remark}

\begin{remark}\label{rem2.7}
It is easy to see that if we demand $f\left(0 \right)=0$ in Theorem \ref{thm2.3} then the series in \eqref{eq2.1} converges uniformly on whole segment $\left[0, \pi \right].$ The same is true for Theorem \ref{thm2.4} if we demand $f\left(\pi \right)=0.$
\end{remark}

\section{Asymptotics of the eigenvalues and the norming constants}
\label{sec3}

As it mentioned above, the case, when $q \in L_{\mathbb{R}}^1 \left[0, \pi \right]$ needs more fine-grained assessment of remainders in the asymptotics of the eigenvalues and the norming constants.

In the papers \cite{Harutyunyan:2016} and \cite{Harutyunyan_Pahlevanyan:2016} the authors found out new asymptotic formulae for the eigenvalues and the norming constants which generalize previously known results. More precisely, the following theorems have been proved:

\begin{theorem}(\cite{Harutyunyan:2016})\label{thm3.1}
Let $q \in L_{\mathbb{R}}^1 \left[0, \pi \right]$ and let $\lambda_n^2\left(q, \alpha, \beta\right)={\mu _n}\left(q, \alpha, \beta\right).$ Then
\begin{enumerate}
\item[(a)] The asymptotic relation $\left(n \to \infty \right)$
\begin{equation}\label{eq3.1}
\lambda_n\left(q, \alpha, \beta\right)=n+\delta_n\left(\alpha, \beta\right)+ \dfrac{[q]}{2\left(n+\delta_n\left(\alpha, \beta\right)\right)}+l_n\left(q, \alpha, \beta\right)+O\left(\dfrac{1}{n^2}\right),
\end{equation}
holds, where $\left[q \right]=\dfrac{1}{\pi}\displaystyle\int_{0}^{\pi}q\left(t\right)dt,$
\begin{equation*}
l_n\left(q, \alpha, \beta\right)=\dfrac{1}{2 \pi \left(n+\delta_n\left(\alpha, \beta\right)\right)} \int_{0}^{\pi} q(x)\cos2\left(n+\delta_n\left(\alpha, \beta\right)\right)xdx, \; \alpha \in \left(0, \pi \right),
\end{equation*}
and
\begin{equation}\label{eq3.2}
l_n=l_n\left(q, \pi, \beta\right)=-\dfrac{1}{2 \pi \left(n+\delta_n\left(\pi, \beta\right)\right)} \int_{0}^{\pi} q(x)\cos2\left(n+\delta_n\left(\pi, \beta\right)\right)xdx.
\end{equation}
The estimate $O\left(\dfrac{1}{n^2}\right)$ of the remainder in \eqref{eq3.1} is uniform in all $\alpha, \beta \in [0, \pi],$ and $q \in {BL}^1_\mathbb{R}\left[0, \pi\right]$ (here and below ${BL}^1_\mathbb{R}\left[0, \pi\right]$ stands for bounded subsets of $L_{\mathbb{R}}^1\left[ {0,\pi } \right]$).
\item[(b)] For $\alpha, \beta \in (0, \pi)$ and for the case $\alpha=\pi,$ $\beta=0$ the function $l,$ defined by the formula
\begin{equation}\label{eq3.3}
l(x)=\displaystyle\sum_{n=2}^{\infty} l_n\left(q, \alpha, \beta\right) \sin \left(n+\delta_n\left(\alpha, \beta\right)\right)x,
\end{equation}
is absolutely continuous on arbitrary segment $\left[a,b\right] \subset \left(0, 2 \pi\right),$ that is $l \in AC\left(0,2\pi\right).$
\end{enumerate}
\end{theorem}

\begin{theorem}(\cite{Harutyunyan_Pahlevanyan:2016})\label{thm3.2}
For norming constants $a_n$ and $b_n$ the following asymptotic formulae hold (when $n \to \infty$):
\begin{multline}\label{eq3.4}
a_n \left( q, \alpha, \beta \right) =\dfrac{\pi}{2} \left[ 1+ \dfrac {2 \, s_n \left( q, \alpha, \beta \right)}{\pi \left[n + \delta \left( \alpha, \beta \right)\right]} + r_n \right] \sin^2 \alpha + \\
+\dfrac{\pi}{2 \left[n + \delta_n(\alpha, \beta)\right]^2} \left[ 1+ \dfrac {2 \, s_n \left( q, \alpha, \beta \right)}{\pi \left[n + \delta \left( \alpha, \beta \right)\right]} + \tilde{r}_n \right]\cos^2 \alpha,
\end{multline}

\begin{multline*}
b_n \left( q, \alpha, \beta \right) =\dfrac{\pi}{2} \left[ 1+ \dfrac {2 \, s_n \left( q, \alpha, \beta \right)}{\pi \left[n + \delta \left( \alpha, \beta \right)\right]} + p_n \right] \sin^2 \beta + \\
+\dfrac{\pi}{2 \left[n + \delta_n(\alpha, \beta)\right]^2} \left[ 1+ \dfrac {2 \, s_n \left( q, \alpha, \beta \right)}{\pi \left[n + \delta \left( \alpha, \beta \right)\right]} + \tilde{p}_n \right]\cos^2 \beta,
\end{multline*}
where
\begin{equation}\label{eq3.5}
s_n = s_n \left( q, \alpha, \beta \right) = -\dfrac {1}{2} \displaystyle \int_{0}^{\pi} \left( \pi - t \right) q\left(t\right) \sin 2 \left[ n + \delta_n\left(\alpha, \beta\right)\right]t dt,
\end{equation}
$r_n = r_n \left( q, \alpha, \beta \right) = O \left( \dfrac {1} {n^2}\right)$ and $\tilde{r}_n=\tilde{r}_n \left( q, \alpha, \beta \right) = O \left( \dfrac {1} {n^2}\right)$ (the same estimate is true for $p_n$ and $\tilde{p}_n$), when $n \to \infty,$ uniformly in $\alpha, \beta \in \left[0, \pi\right]$ and $q \in {BL}^1_\mathbb{R}\left[0, \pi\right].$
\end{theorem}

\begin{theorem}(\cite{Harutyunyan_Pahlevanyan:2016})\label{thm3.3}
For both $\alpha, \beta \in (0, \pi)$ and $\alpha = \pi, \; \beta = 0$ cases the function $s,$ defined as the series
\begin{equation}\label{eq3.6}
s \left(x \right) = \displaystyle \sum_{n=2}^{\infty} \dfrac {s_n}{n + \delta_n \left( \alpha, \beta \right)} \cos \left[n + \delta_n \left( \alpha, \beta \right)\right]x
\end{equation}
is absolutely continuous function on arbitrary segment $\left[a,b\right] \subset \left(0, 2 \pi\right),$ i.e. $s \in AC \left(0, 2 \pi\right).$
\end{theorem}

The proofs of the Theorem \ref{thm3.1} and Theorem \ref{thm3.3} do not cover the case $\alpha=\pi, \; \beta \in \left(0,\pi \right).$ In the forthcoming theorems using Theorem \ref{thm2.3} and Theorem \ref{thm2.4} we handle this case as well.
\begin{theorem}\label{thm3.4}
The function $l,$ defined by the formula
\begin{equation*}
l \left(x, \beta \right)=\displaystyle\sum_{n=2}^{\infty} l_n\left(q, \pi, \beta\right) \sin \left(n+\delta_n\left(\pi, \beta\right)\right)x
\end{equation*}
is absolutely continuous on arbitrary segment $\left[a,b\right] \subset \left(0, 2 \pi\right),$ i.e. $l \in AC \left(0, 2\pi \right).$
\end{theorem}
\begin{proof}
Denote $\sigma\left(x\right)=\displaystyle\int_{0}^{x} q\left( t \right)dt$ and write $l_n\left(q, \pi, \beta \right)$ in the following form:
\begin{multline}\label{eq3.7}
l_n\left(q, \pi, \beta\right)=-\dfrac{1}{2 \pi \left(n+\delta_n\left(\pi, \beta\right)\right)} \int_{0}^{\pi} q(x)\cos2\left(n+\delta_n\left(\pi, \beta\right)\right)xdx=\\
=-\dfrac{1}{2 \pi \left(n+\delta_n\left(\pi, \beta\right)\right)} \int_{0}^{\pi}\cos2\left(n+\delta_n\left(\pi, \beta\right)\right)xd\sigma\left(x\right)=-\dfrac{\sigma \left(\pi\right)\cos 2\pi \delta_n\left(\pi, \beta\right)}{2 \pi \left(n+\delta_n\left(\pi, \beta\right)\right)}-\\
-\dfrac{1}{\pi}\int_{0}^{\pi}\sigma \left(x\right)\sin2\left(n+\delta_n\left(\pi, \beta\right)\right)xdx=-\dfrac{\sigma \left(\pi\right)\cos 2\pi \delta_n\left(\pi, \beta\right)}{2 \pi \left(n+\delta_n\left(\pi, \beta\right)\right)}-\\
-\dfrac{1}{2\pi}\int_{0}^{2\pi}\sigma_1 \left(x\right)\sin\left(n+\delta_n\left(\pi, \beta\right)\right)xdx,
\end{multline}
where $\sigma_1 \left( x \right)\equiv \sigma \left( \dfrac{x}{2} \right)$ is the absolutely continuous function on $\left[0, 2\pi\right].$

It is easy to see from \eqref{eq1.10} (for the details see \cite{Harutyunyan:2008}), that for $\beta \in \left(0,\pi\right)$ we have
\begin{equation}\label{eq3.8}
\delta_n \left(\pi, \beta \right)=\dfrac{1}{2}+\dfrac{\cot\beta}{\pi\left(n+\frac{1}{2}\right)}+ O\left(\dfrac{1}{n^2}\right)\,{\cot\beta} =\dfrac{1}{2}+O\left(\dfrac{1}{n}\right),
\end{equation}
and consequently,
\begin{equation}\label{eq3.9}
\cos 2 \pi \delta_n\left(\pi, \beta\right)=-1+d_n, \; \sin 2 \pi \delta_n\left(\pi, \beta\right)=e_n,
\end{equation}
where $d_n=O\left(\dfrac{1}{n^2}\right), \; e_n=O\left(\dfrac{1}{n}\right).$

Therefore, $l\left(x, \beta \right)$ can be represented as a sum of three functions
\begin{equation}\label{eq3.10}
l\left(x, \beta\right)=l_1\left(x, \beta\right)+l_2\left(x, \beta\right)+l_3\left(x, \beta\right),
\end{equation}
where
\begin{equation}\label{eq3.11}
l_1\left(x, \beta\right)=\dfrac{\sigma\left(\pi\right)}{2 \pi}\displaystyle\sum_{n=2}^{\infty} \dfrac{\sin \left(n+\delta_n\left(\pi, \beta\right)\right)x}{\left(n+\delta_n\left(\pi, \beta\right)\right)},
\end{equation}
\begin{equation}\label{eq3.12}
l_2\left(x, \beta\right)=-\dfrac{\sigma\left(\pi\right)}{2 \pi}\displaystyle\sum_{n=2}^{\infty} {d_n}\dfrac{\sin \left(n+\delta_n\left(\pi, \beta\right)\right)x}{\left(n+\delta_n\left(\pi, \beta\right)\right)},
\end{equation}
\begin{equation}\label{eq3.13}
l_3\left(x, \beta\right)=-\dfrac{1}{2 \pi} \displaystyle\sum_{n=2}^{\infty} f_n \sin \left(n+\delta_n\left(\pi, \beta\right)\right)x,
\end{equation}
and $f_n=\displaystyle\int_{0}^{2\pi} \sigma_1 \left(t\right)\sin\left(n+\delta_n\left(\pi, \beta\right)\right)tdt.$

Since $f_n=\displaystyle\int_{0}^{\pi} \sigma_1 \left(t\right)\sin\left(n+\delta_n\left(\pi, \beta\right)\right)tdt+\displaystyle\int_{\pi}^{2\pi} \sigma_1 \left(t\right)\sin\left(n+\delta_n\left(\pi, \beta\right)\right)tdt$ and \\
\begin{multline*}
\int_{\pi}^{2\pi} \sigma_1 \left(t\right)\sin\left(n+\delta_n\left(\pi, \beta\right)\right)tdt=
\int_{-2\pi}^{-\pi} -\sigma_1 \left(-t\right)\sin\left(n+\delta_n\left(\pi, \beta\right)\right)tdt=\\
=\int_{0}^{\pi} -\sigma_1 \left(2\pi-t\right)\sin\left(n+\delta_n\left(\pi, \beta\right)\right)\left(t-2\pi\right)dt=\\
=\int_{0}^{\pi} \sigma_1 \left(2\pi-t\right)\left(\left(1-d_n\right)\sin\left(n+\delta_n\left(\pi, \beta\right)\right)t+e_n \cos\left(n+\delta_n\left(\pi, \beta\right)\right)t \right)dt=\\
=\int_{0}^{\pi} \sigma \left(\pi-\dfrac{t}{2}\right)\left(\left(1-d_n\right)\sin\left(n+\delta_n\left(\pi, \beta\right)\right)t+e_n \cos\left(n+\delta_n\left(\pi, \beta\right)\right)t \right)dt,
\end{multline*}
then
\begin{multline}\label{eq3.14}
f_n=\int_{0}^{\pi} \left(\sigma \left(\dfrac{t}{2}\right)+\sigma \left(\pi-\dfrac{t}{2}\right)\right)\sin\left(n+\delta_n\left(\pi, \beta\right)\right)tdt-\\
-d_n \int_{0}^{\pi} \sigma \left(\pi-\dfrac{t}{2}\right)\sin\left(n+\delta_n\left(\pi, \beta\right)\right)tdt+\\
+e_n \int_{0}^{\pi} \sigma \left(\pi-\dfrac{t}{2}\right)\cos\left(n+\delta_n\left(\pi, \beta\right)\right)t.
\end{multline}
It noteworthy, that $\delta_n\left(\alpha, \beta\right)$ are defined only for $n \geq 2,$ that is why, we will write $\lambda_0 \left(0,\pi,\beta\right),$ $\lambda_1 \left(0,\pi,\beta\right)$ and $\lambda_n \left(0,\pi,\beta\right)=n+\delta_n\left(\pi, \beta\right)$ for all $n \geq 2.$
Taking into account that the system of functions
$$\left\{\varphi_{n} \left(x \right)\right\}_{n=0}^\infty=\left\{\dfrac{\sin \lambda_n \left(0,\pi,\beta\right)x}{\lambda_n \left(0,\pi,\beta\right)}\right\}_{n=0}^\infty = \left\{\dfrac{\sin \lambda_n \left(0,\pi,\beta\right)x}{\lambda_n \left(0,\pi,\beta\right)}\right\}_{n=0}^1 \cup \left\{\dfrac{\sin \left(n+\delta_n\left(\pi, \beta \right) \right)x}{n+\delta_n\left(\pi, \beta \right)}\right\}_{n=2}^\infty$$
are the eigenfunctions of the problem $L\left(0, \pi, \beta\right)$ and using Theorem \ref{thm2.3}, we get
\begin{multline}\label{eq3.15}
\sigma \left(\dfrac{x}{2}\right)+\sigma \left(\pi-\dfrac{x}{2}\right) =\sigma_{2}\left(x\right)+\\
+\sum\limits_{n = 2}^\infty \dfrac{\displaystyle\int_{0}^{\pi}\left(\sigma \left(\dfrac{t}{2}\right)+\sigma \left(\pi-\dfrac{t}{2}\right)\right)\sin \left(n+\delta_n\left(\pi, \beta \right) \right)tdt}{\displaystyle\int_{0}^{\pi}\sin^{2}\left(n+\delta_n\left(\pi, \beta \right) \right) tdt} \sin \left(n+\delta_n\left(\pi, \beta \right) \right)x
\end{multline}
where the series converges uniformly on arbitrary segment $\left[a, \pi \right] \subset \left(0, \pi \right]$ and
\begin{equation*}
\sigma_{2}\left(x\right):= \sum\limits_{n = 0}^1 \dfrac{\displaystyle\int_{0}^{\pi}\left(\sigma \left(\dfrac{t}{2}\right)+\sigma \left(\pi-\dfrac{t}{2}\right)\right)\sin \lambda_n \left(0,\pi,\beta\right)tdt}{\displaystyle\int_{0}^{\pi}\sin^{2}\lambda_n \left(0,\pi,\beta\right)tdt}\sin \lambda_n \left(0,\pi,\beta\right)x
\end{equation*}
Using \eqref{eq3.8} and \eqref{eq3.9}, we calculate
\begin{multline}\label{eq3.16}
\int\limits_0^\pi {\sin^2 \left(n+{\delta_n}\left( \pi, \beta \right)\right)t}dt=\\
=\dfrac{\pi}{2}-\dfrac{\sin 2 \pi \left(n+\delta_n (\pi, \beta)\right)}{4\left(n+\delta_n (\pi, \beta)\right)}=\dfrac{\pi}{2} - \dfrac{e_n}{4\left(n+\delta_n (\pi, \beta)\right)}.
\end{multline}
From \eqref{eq3.16}, it is easy to see that
$$\dfrac{1}{\int\limits_0^\pi {\sin^2 \left(n+{\delta_n}\left( \pi, \beta \right) \right)t}dt} = \dfrac{2}{\pi}+{g_n},$$
where $g_n=\dfrac{2 e_n}{\pi\left(2\pi\left(n+\delta_n\left(\pi, \beta \right) \right)-e_n \right)}=O\left(\dfrac{1}{n^2}\right).$

Now we can write \eqref{eq3.15} in the form
\begin{multline}\label{eq3.17}
\sum\limits_{n = 2}^\infty \dfrac{2}{\pi}\displaystyle\int_{0}^{\pi}\left(\sigma \left(\dfrac{t}{2}\right)+\sigma \left(\pi-\dfrac{t}{2}\right)\right)\sin \left(n+\delta_n\left(\pi, \beta \right) \right)tdt \sin \left(n+\delta_n\left(\pi, \beta \right) \right)x=\\
=-\sum\limits_{n = 2}^\infty {g_n}\displaystyle\int_{0}^{\pi}\left(\sigma \left(\dfrac{t}{2}\right)+\sigma \left(\pi-\dfrac{t}{2}\right)\right)\sin \left(n+\delta_n\left(\pi, \beta \right) \right)tdt \sin \left(n+\delta_n\left(\pi, \beta \right) \right)x+\\
+\sigma \left(\dfrac{x}{2}\right)+\sigma \left(\pi-\dfrac{x}{2}\right) -\sigma_2\left(x\right),
\end{multline}
where the series converge uniformly on arbitrary segment $\left[a, \pi \right] \subset \left(0, \pi \right].$

Hence, from \eqref{eq3.13}, \eqref{eq3.14}, \eqref{eq3.17} we obtain that for arbitrary $x \in \left(0, \pi \right]$
\begin{multline}\label{eq3.18}
l_3\left(x, \beta\right) = \dfrac{1}{4}\left(-\sigma \left(\dfrac{x}{2}\right)-\sigma \left(\pi-\dfrac{x}{2}\right)+\sigma_2\left(x\right)\right)+\\
+\displaystyle\sum_{n=2}^{\infty} \dfrac{g_n}{4} \displaystyle\int_{0}^{\pi}\left(\sigma \left(\dfrac{t}{2}\right)+\sigma \left(\pi-\dfrac{t}{2}\right)\right)\sin \left(n+\delta_n\left(\pi, \beta \right) \right)tdt \sin \left(n+\delta_n\left(\pi, \beta \right) \right)x+\\
+\dfrac{1}{2\pi}\displaystyle\sum_{n=2}^{\infty} {d_n} \int_{0}^{\pi} \sigma \left(\pi-\dfrac{t}{2}\right)\sin\left(n+\delta_n\left(\pi, \beta\right)\right)tdt \sin \left(n+\delta_n\left(\pi, \beta\right)\right)x - \\
-\dfrac{1}{2\pi}\displaystyle\sum_{n=2}^{\infty} {e_n} \int_{0}^{\pi} \sigma \left(\pi-\dfrac{t}{2}\right)\cos\left(n+\delta_n\left(\pi, \beta\right)\right)tdt \sin \left(n+\delta_n\left(\pi, \beta\right)\right)x.
\end{multline}
Since $d_n=O\left(\dfrac{1}{n^2}\right),$ $e_n \displaystyle\int_{0}^{\pi} \sigma \left(\pi-\dfrac{t}{2}\right)\cos\left(n+\delta_n\left(\pi, \beta\right)\right)tdt=O\left(\dfrac{1}{n^2}\right),$ $g_n=O\left(\dfrac{1}{n^2}\right),$ then $l_3 \in AC\left(0, \pi\right].$ But on the other hand, since (see \eqref{eq3.13} and \eqref{eq3.9})
\begin{multline*}
l_3\left(2\pi-x, \beta\right)=l_3\left(x, \beta\right)+\dfrac{1}{2 \pi} \displaystyle\sum_{n=2}^{\infty} d_n f_n \sin \left(n+\delta_n\left(\pi, \beta\right)\right)x-\\
-\dfrac{1}{2 \pi} \displaystyle\sum_{n=2}^{\infty} e_n f_n \cos \left(n+\delta_n\left(\pi, \beta\right)\right)x,
\end{multline*}
then $l_3 \in AC\left[\pi, 2\pi\right)$ and consequently $l_3 \in AC\left(0, 2\pi\right).$ \\
Since $\displaystyle\sum_{n=2}^{\infty} \dfrac{\sin \left(n+\delta_n\left(\pi, \beta\right)\right)x}{\left(n+\delta_n\left(\pi, \beta\right)\right)}$ is absolutely continuous function on $\left(0, 2 \pi\right)$ (see section \ref{app}, Lemma \ref{lem7.1}), then $l_1 \in AC\left(0, 2\pi\right).$ \\
Since $d_n=O\left(\dfrac{1}{n^2}\right),$ then the series in \eqref{eq3.12} and its first derivative converges absolutely and uniformly on $\left[0, 2\pi\right]$ and therefore $l_2 \in AC\left[0, 2\pi\right].$ This completes the proof.
\end{proof}

A similar assertion is true for remainders of the norming constants. More precisely, we prove the following theorem:
\begin{theorem}\label{thm3.5}
The function $s,$ defined by the formula
\begin{equation*}
s \left(x, \beta \right) = \displaystyle \sum_{n=2}^{\infty} \dfrac {s_n\left(q, \pi, \beta \right)}{n + \delta_n \left( \pi, \beta \right)} \cos \left(n + \delta_n \left( \pi, \beta \right)\right)x
\end{equation*}
is absolutely continuous on $\left[0, 2 \pi\right],$ i.e. $s \in AC \left[0, 2 \pi\right].$
\end{theorem}
\begin{proof}
Since (see \eqref{eq3.8})
\begin{multline*}
\sin 2\left(n+\delta_n \left(\pi,\beta\right)\right)x-\sin 2\left(n+\dfrac{1}{2}\right)x=\\
=-2\cos\left(n+\delta_n \left(\pi,\beta\right)+\dfrac{1}{2}\right)x \sin\left(\delta_n \left(\pi,\beta\right)-\dfrac{1}{2}\right)x = O\left(\dfrac{1}{n}\right),
\end{multline*}
\begin{multline*}
\cos \left(n+\delta_n \left(\pi,\beta\right)\right)x-\cos\left(n+\dfrac{1}{2}\right)x=\\
=-2\sin\left(n+\dfrac{2\delta_n \left(\pi,\beta\right)+1}{4}\right)x \sin\left(\dfrac{2\delta_n \left(\pi,\beta\right)-1}{4}\right)x = O\left(\dfrac{1}{n}\right)
\end{multline*}
uniformly on $\left[0,2\pi\right],$ then $s \in AC\left[0,2\pi\right]$ is equivalent to $\tilde{s} \in AC\left[0,2\pi\right],$ where $\tilde{s}$ defined by the formula
\begin{equation}\label{eq3.19}
\tilde{s}(x,\beta) := \displaystyle \sum_{n=0}^{\infty} \dfrac {\tilde{s_n}\left(q,\pi,\beta\right)}{n + \frac{1}{2}} \cos \left(n + \dfrac{1}{2}\right)x,
\end{equation}
where
\begin{equation*}
\tilde{s_n} \left( q, \pi, \beta \right) = -\dfrac {1}{2} \displaystyle \int_{0}^{\pi} \left( \pi - t \right) q\left(t\right) \sin 2 \left( n + \dfrac{1}{2}\right)t dt.
\end{equation*}
Denote ${\sigma_3}\left(x\right)=\displaystyle\int_{0}^{x} \left(\pi-t\right)q\left( t \right)dt$ and write $\dfrac {\tilde{s_n}\left(q,\pi,\beta\right)}{n + \frac{1}{2}}$ in the following form:
\begin{multline}\label{eq3.20}
\dfrac {\tilde{s_n}\left(q,\pi,\beta\right)}{n + \frac{1}{2}}=-\dfrac{1}{2\left(n + \frac{1}{2}\right)} \displaystyle \int_{0}^{\pi} \left( \pi - t \right) q\left(t\right) \sin 2 \left( n + \dfrac{1}{2}\right)tdt=\\
=-\dfrac{1}{2\left(n + \frac{1}{2}\right)}\displaystyle \int_{0}^{\pi} \sin 2 \left( n + \dfrac{1}{2}\right)td{\sigma_3}\left(t\right)=\displaystyle \int_{0}^{\pi}{\sigma_3}\left(t\right)\cos 2 \left( n + \dfrac{1}{2}\right)tdt=\\
=\dfrac{1}{2}\displaystyle \int_{0}^{2\pi}{\sigma_4}\left(t\right)\cos \left( n + \dfrac{1}{2}\right)tdt,
\end{multline}
where ${\sigma_4}\left(x\right)\equiv{\sigma_3}\left(\dfrac{x}{2}\right)$ is the absolutely continuous function on $\left[0, 2\pi\right].$

Since
\begin{multline*}
\displaystyle\int_{\pi}^{2\pi}{\sigma_4}\left(t\right)\cos\left(n+\dfrac{1}{2}\right)tdt=
\displaystyle\int_{-2\pi}^{-\pi}{\sigma_4}\left(-t\right)\cos\left(n+\dfrac{1}{2}\right)tdt=\\
=\displaystyle\int_{0}^{\pi}{\sigma_4}\left(2\pi-t\right)\cos\left(n+\dfrac{1}{2}\right)
\left(t-2\pi\right)dt=-\displaystyle\int_{0}^{\pi}{\sigma_4}\left(2\pi-t\right)\cos\left(n+\dfrac{1}{2}\right)tdt,
\end{multline*}
then the integral in \eqref{eq3.20} can be written in the following form
\begin{multline}\label{eq3.21}
\dfrac{1}{2}\displaystyle\int_{0}^{2\pi}{\sigma_4}\left(t\right)\cos\left(n+\dfrac{1}{2}\right)tdt=\\
=\dfrac{1}{2}\displaystyle\int_{0}^{\pi}\left({\sigma_4}\left(t\right)-{\sigma_4}\left(2\pi-t\right)\right)
\cos\left(n+\dfrac{1}{2}\right)tdt.
\end{multline}
Taking into account that the system of functions $\left\{\cos \left(n+\dfrac{1}{2}\right)x\right\}_{n=0}^\infty$
are the eigenfunctions of the problem $L\left(0, \dfrac{\pi}{2}, 0\right)$ and using Theorem \ref{thm2.4} and Remark \ref{rem2.7}, we receive
\begin{multline}\label{eq3.22}
{\sigma_4}\left(x\right)-{\sigma_4}\left(2\pi-x\right)=\\
=\dfrac{2}{\pi}\displaystyle \sum_{n=0}^{\infty} \displaystyle \int_{0}^{\pi} \left( {\sigma_4}\left(t\right)-{\sigma_4}\left(2\pi-t\right) \right) \cos \left(n + \dfrac{1}{2}\right)tdt \cos \left(n + \dfrac{1}{2}\right)x,
\end{multline}
where the series converges uniformly on $\left[0,\pi\right].$

From \eqref{eq3.19}, \eqref{eq3.20}, \eqref{eq3.21} and \eqref{eq3.22} we get that
\begin{equation*}
\tilde{s}\left(x\right)=\dfrac{\pi}{4}\left({\sigma_4}\left(x\right)-{\sigma_4}\left(2\pi-x\right)\right), \; x \in \left[0, \pi \right],
\end{equation*}
but on the other hand (see \eqref{eq3.19})
\begin{equation*}
\tilde{s} \left(2\pi-x, \beta \right)=-\tilde{s} \left(x, \beta \right)
\end{equation*}
and therefore $\tilde{s}$ is an absolutely continuous function on $\left[0,2\pi\right].$ This completes the proof.
\end{proof}
\begin{remark}\label{rem3.6}
We have proved that in the direct problem ($q$ and $\beta$ are given) the properties \eqref{eq1.11}, \eqref{eq1.13} of the eigenvalues and the properties \eqref{eq1.12}, \eqref{eq1.14} of the norming constants take place.
\end{remark}
At the end of this section we formulate two theorems on more precise asymptotic formulae for eigenvalues and norming constants, when the potential $q$ is an absolutely continuous function.

\begin{theorem}\label{thm3.7}
Let $q \in AC_{\mathbb{R}} \left[0, \pi \right].$ Then the asymptotic relation $\left(n \to \infty \right)$ for eigenvalues
\begin{equation}\label{eq3.23}
\lambda_{n} \left(q, \pi, \beta \right)=n+\delta_n \left(\pi, \beta \right)+\dfrac {\left[q\right]}{2 \left(n+\delta_n \left(\pi, \beta \right) \right)}+l_{n} \left(q, \pi, \beta \right)+O\left(\dfrac{1}{n^3}\right),
\end{equation}
holds, where $\left[q \right]=\dfrac{1}{\pi}\displaystyle\int_{0}^{\pi}q\left(t\right)dt,$
\begin{equation}\label{eq3.24}
l_n=l_n\left(q, \pi, \beta\right)=\dfrac{1}{4 \pi \left(n+\delta_n\left(\pi, \beta\right)\right)^2} \int_{0}^{\pi} q'\left(x\right)\sin 2\left(n+\delta_n\left(\pi, \beta\right)\right)xdx.
\end{equation}
The estimate $O\left(\dfrac{1}{n^3}\right)$ of the remainder in \eqref{eq3.23} is uniform in all $\beta \in [0, \pi]$ and $q, q' \in {BL}^1_\mathbb{R}\left[0, \pi\right],$ and the functions
\begin{equation}\label{eq3.25}
l\left(x, \beta \right)=\displaystyle\sum_{n=2}^{\infty} l_n\left(q, \pi, \beta\right) \sin \left(n+\delta_n\left(\pi, \beta\right)\right)x,
\end{equation}
\begin{equation}\label{eq3.26}
l'\left(x, \beta \right)=\displaystyle\sum_{n=2}^{\infty} l_n\left(q, \pi, \beta\right) \left(n+\delta_n\left(\pi, \beta\right)\right) \cos \left(n+\delta_n\left(\pi, \beta\right)\right)x
\end{equation}
are absolutely continuous on arbitrary segment $\left[a, b\right] \subset \left(0, 2 \pi\right),$ that is $l, l' \in AC\left(0,2\pi\right).$
\end{theorem}

\begin{theorem}\label{thm3.8}
Let $q \in AC_{\mathbb{R}} \left[0, \pi \right].$ Then the asymptotic relation $\left(n \to \infty \right)$ for norming constants
\begin{equation}\label{eq3.27}
a_{n} \left(q, \pi, \beta \right)=\dfrac {\pi}{2 \left(n+{\delta_n} \left(\pi, \beta \right) \right)^2} \left(1+\dfrac{\left[q \right]_{\beta}}{2 \left(n+\delta_n \left(\pi, \beta \right) \right)^2} +\dfrac{2 \, s_{n}}{\pi \left(n+\delta_n \left(\pi, \beta \right) \right)^2}+O\left(\dfrac{1}{n^3}\right)\right),
\end{equation}
holds, where $\left[q \right]_{\beta}=\dfrac{5}{\pi}\displaystyle\int_{0}^{\pi}q\left(t\right)dt+2\left( q\left(0 \right)+\cot \beta \right),$
\begin{equation}\label{eq3.28}
s_n=s_n\left(q, \pi, \beta\right)=\dfrac{1}{4} \int_{0}^{\pi} \left(\pi-t\right) q'\left(t\right)\cos 2\left(n+\delta_n\left(\pi, \beta\right)\right)tdt.
\end{equation}
The estimate $O\left(\dfrac{1}{n^3}\right)$ of the remainder in \eqref{eq3.27} is uniform in all $\beta \in [0, \pi]$ and $q, q' \in {BL}^1_\mathbb{R}\left[0, \pi\right],$ and the functions
\begin{equation}\label{eq3.29}
s\left(x, \beta \right)=\displaystyle\sum_{n=2}^{\infty} \dfrac{s_n\left(q, \pi, \beta\right)}{\left(n+\delta_n\left(\pi, \beta\right)\right)^2} \cos \left(n+\delta_n\left(\pi, \beta\right)\right)x,
\end{equation}
\begin{equation}\label{eq3.30}
s'\left(x, \beta \right)=-\displaystyle\sum_{n=2}^{\infty} \dfrac{s_n\left(q, \pi, \beta\right)}{n+\delta_n\left(\pi, \beta\right)} \sin \left(n+\delta_n\left(\pi, \beta\right)\right)x
\end{equation}
are absolutely continuous on arbitrary segment $\left[a, b\right] \subset \left(0, 2 \pi\right),$ that is $s, s' \in AC\left(0,2\pi\right).$
\end{theorem}

The proof of these theorems is based on the fact that for the case $q \in AC_{\mathbb{R}} \left[0, \pi \right]$ more precise asymtotic formula $\left(n \to \infty \right)$ for the eigenfunctions
\begin{multline*}
\varphi \left(x, \mu_{n} \left(q, \pi, \beta \right), \pi \right)=\dfrac{\sin \lambda_{n} x}{\lambda_{n}}-\dfrac{\cos \lambda_{n} x}{2 \lambda_{n}^2}\displaystyle\int_{0}^{x} q \left(t \right)dt+\\
+\dfrac{\sin \lambda_{n} x}{4 \lambda_{n}^3} \left(q\left(x \right)+q\left(0 \right)-\dfrac{1}{2}\left(\displaystyle\int_{0}^{x} q \left(t \right)dt \right)^2\right)+\\
+\dfrac{1}{4 \lambda_{n}^3}\displaystyle\int_{0}^{x} q' \left(t \right) \sin \lambda_{n} \left(x-2t \right) dt+O\left(\dfrac{1}{n^4}\right).
\end{multline*}
holds (can be obtained using integration by parts method). The rest part of the proof is based on the same technique as for general case $q \in L_{\mathbb{R}}^1 \left[0, \pi \right]$ (see \cite[the proof of theorem \ref{thm3.1}]{Harutyunyan:2016} and \cite[the proofs of theorems \ref{thm3.2}--\ref{thm3.3}]{Harutyunyan_Pahlevanyan:2016}), that is why we omit it.

\section{Derivation of an analogue of the Gelfand-Levitan equation}
\label{sec4}

Transformation (transmutation) operators are played an important role in the theory of the inverse Sturm-Liouville problems.

The following assertion can be found in \cite{Freiling_Yurko:2001} (see, also \cite{Povzner:1948,Marchenko:1977}):

\begin{theorem}(\cite{Freiling_Yurko:2001})\label{thm4.1}
For the function $\varphi_{\pi} \left(x, \mu \right)$ the following representation holds
\begin{equation}\label{eq4.1}
\varphi_{\pi} \left(x, \mu \right) =\left(\mathbb{I}+\mathbb{P}\right)\dfrac{\sin \lambda x}{\lambda} := \dfrac{\sin \lambda x}{\lambda}+{\int_{0}^{x}}P(x,t)\dfrac{\sin \lambda t} {\lambda}dt,
\end{equation}
where $P\left(x,t\right),$ $0 \leq t \leq x \leq \pi,$ is a real continuous function with the same smoothness as $\displaystyle{\int_{0}^{x}}q(t)dt$, and
\begin{equation}\label{eq4.2}
P(x,x)=\dfrac{1}{2}{\int_{0}^{x}}q(t)dt, \; P(x,0)=0.
\end{equation}
\end{theorem}
The proof of the theorem \ref{thm4.1} was given for the case $q \in L_{\mathbb{R}}^2 \left[0, \pi \right],$ but it can be easily done for the case $q \in L_{\mathbb{R}}^1 \left[0, \pi \right],$ without any changes.

Further, in this section, our goal is to derive the analogue of Gelfand-Levitan equation for our case.
\begin{lemma}\label{lem4.2}
Let us given the numbers $\left\{\mu_n\right\}_{n \geq 0}$ and $\left\{a_n \right\}_{n \geq 0},$ which satisfy the conditions of the Theorem \ref{thm1.1}, such that $\lambda_n \neq 0,$ $n=0,1,2,\dots$ and $\lambda_n \left(0, \pi, \beta \right) \neq 0,$ $n=0,1,2,\dots.$ Then the function $H,$ defined by the formula
\begin{equation}\label{eq4.3}
H \left(t \right)=\sum\limits_{n = 0}^{\infty}\left(\dfrac{1}{a_n} \, \dfrac{\cos \lambda_{n} t}{\lambda_{n}^2}- \dfrac{1}{{a_n} \left( 0, \pi, \beta \right)} \, \dfrac{\cos \lambda_{n} \left(0, \pi, \beta \right) t} {\lambda_{n}^{2} \left(0, \pi, \beta \right)}\right),
\end{equation}
is absolutely continuous on arbitrary segment $\left[a, b \right] \subset \left(0, 2 \pi \right)$ $\left( H\left(\cdot\right) \in AC(0, 2\pi) \right).$
\end{lemma}
\begin{proof}
Denote
\begin{equation*}
k_n :={a_n}\lambda_n^2, \; k_n \left(0, \pi, \beta \right):={{a_n} \left( 0, \pi, \beta \right)}\left(n+{\delta_n}(\pi, \beta)\right)^2, \; \epsilon_n:=\dfrac{c}{2(n+{\delta_n}(\pi, \beta))}+l_n.
\end{equation*}
Since $\lambda_{n} \left(0, \pi, \beta \right)=n+\delta_{n} \left(\pi, \beta \right),$ $n=2,3,\dots,$ then, we can write down the general term of the series \eqref{eq4.3} in the following form:
\begin{multline}\label{eq4.4}
\dfrac{1}{a_n} \, \dfrac{\cos \lambda_{n} t}{\lambda_{n}^2}- \dfrac{1}{{a_n} \left( 0, \pi, \beta \right)} \, \dfrac{\cos \left(n+\delta_n (\pi, \beta) \right)t}{{\left(n+\delta_n (\pi, \beta)\right)}^2}=\\
=\dfrac{1}{k_n}\left( \cos \left(n+{\delta_n}\left(\pi, \beta\right)+\epsilon_n\right)t-\cos \left(n+{\delta_n}\left(\pi, \beta\right)\right)t \right)+\\
+\left(\dfrac{1}{k_n}-\dfrac{1}{k_n \left(0, \pi, \beta \right)}\right)\cos \left(n+{\delta_n}(\pi, \beta)\right)t.
\end{multline}
Using \eqref{eq1.11}, \eqref{eq1.12} and the condition $l_n=o\left(\dfrac{1}{n}\right)$ we can write $\dfrac{1}{k_n}$ in the form:
\begin{multline}\label{eq4.5}
\dfrac{1}{k_n}=\dfrac{1}{a_n \lambda_n^2}=\dfrac{2}{\pi \left( 1+ \dfrac {2 \, s_n}{\pi \left(n + \delta \left( \pi, \beta \right)\right)}\right) \left(1+\dfrac{c+r_n}{\left(n+\delta_n\left(\pi, \beta\right)\right)^2} \right)}=\\
=\dfrac{2}{\pi}-\dfrac{4}{\pi^2} \, \dfrac{s_n}{n+\delta_n \left(\pi, \beta \right)}+h_n,
\end{multline}
where $r_n = o\left(1 \right),$ $h_n=O\left(\dfrac{1}{n^2}\right).$ On the other hand,
\begin{equation*}
\dfrac{1}{k_n \left(0,\pi,\beta \right)}= \dfrac{1}{{{a_n} \left( 0, \pi, \beta \right)}\left(n+{\delta_n}(\pi, \beta)\right)^2}=\dfrac{2}{\pi}+g_n, \;\; g_n = \left(\dfrac{1}{n^2}\right).
\end{equation*}
Therefore
\begin{equation}\label{eq4.6}
\dfrac{1}{k_n}-\dfrac{1}{k_n \left(0, \pi, \beta \right)}=-\dfrac{4}{\pi^2} \, \dfrac{s_n}{n+\delta_n \left(\pi, \beta \right)}+h_n-g_n.
\end{equation}
And finally, using elementary trigonometric identities and Maclaurin expansions of $\cos$ and $\sin$ functions around the point $0,$ we get that
\begin{multline}\label{eq4.7}
\cos \left(n+{\delta_n}\left(\pi, \beta\right)+\epsilon_n\right)t-\cos \left(n+{\delta_n}\left(\pi, \beta\right)\right)t=\\
=-\epsilon_n t \sin \left(n+\delta_n \left(\pi, \beta \right) \right)t + m_n \left(t \right) \cos \left(n+\delta_n \left(\pi, \beta \right) \right)t - \\
-z_n \left(t \right) \sin \left(n+\delta_n \left(\pi, \beta \right) \right)t,
\end{multline}
where $m_n \left(t \right)=O\left(\dfrac{1}{n^2}\right)$ and $z_n \left(t \right)=O\left(\dfrac{1}{n^3}\right)$ uniformly with respect to $t \in \left[0, 2 \pi \right].$

Now by taking into account \eqref{eq4.4}, \eqref{eq4.5}, \eqref{eq4.6} and \eqref{eq4.7} we can expand $H\left(t \right)$ in the following form:
\begin{equation}\label{eq4.8}
H\left(t \right)=\displaystyle \sum\limits_{n = 0}^{1}\left(\dfrac{1}{a_n} \dfrac{\cos \lambda_{n} t}{\lambda_{n}^2}- \dfrac{1}{{a_n} \left( 0, \pi, \beta \right)} \dfrac{\cos \lambda_{n}\left(0, \pi, \beta\right) t} {\lambda_{n}^{2}\left(0, \pi, \beta\right)}\right)+\sum\limits_{i=1}^{5}{H_i}\left(t\right),
\end{equation}
where
\begin{equation}\label{eq4.9}
H_1 \left(t \right)=-\dfrac{ct}{\pi} \displaystyle \sum\limits_{n = 2}^\infty {\dfrac{\sin \left( n+{\delta_n}\left(\pi, \beta\right)\right)t}{n+{\delta_n}\left(\pi, \beta\right)}},
\end{equation}
\begin{equation}\label{eq4.10}
H_2 \left(t  \right)=-\dfrac{2t}{\pi} \displaystyle \sum\limits_{n = 2}^\infty {{l_n} {\sin \left( n+{\delta_n}\left(\pi, \beta \right)\right)t}}=-\dfrac{2t}{\pi} \, l\left(t \right),
\end{equation}
\begin{multline}\label{eq4.11}
H_3 \left(t \right)= -\displaystyle \sum_{n=2}^{\infty} \left(-\dfrac{4}{\pi^2} \, \dfrac {s_n}{n + \delta_n \left( \pi, \beta \right)}+ h_n\right) \epsilon_n t \sin \left(n + \delta_n \left( \pi, \beta \right)\right)t + \\ \displaystyle \sum_{n=2}^{\infty} \left(\dfrac{2}{\pi}-\dfrac{4}{\pi^2} \, \dfrac {s_n}{n + \delta_n \left( \pi, \beta \right)}+ h_n\right) m_n \left(t \right) \cos \left(n + \delta_n \left( \pi, \beta \right)\right)t - \\
-\displaystyle \sum_{n=2}^{\infty} \left(\dfrac{2}{\pi}-\dfrac{4}{\pi^2} \, \dfrac {s_n}{n + \delta_n \left( \pi, \beta \right)}+ h_n\right) z_n \left(t \right) \sin \left(n + \delta_n \left( \pi, \beta \right)\right)t,
\end{multline}
\begin{equation}\label{eq4.12}
H_4 \left(t \right)=-\dfrac{4}{\pi^2} \displaystyle \sum_{n=2}^{\infty} \dfrac {s_n}{n + \delta_n \left( \pi, \beta \right)} \cos \left(n + \delta_n \left( \pi, \beta \right)\right)t=-\dfrac{4}{\pi^2} \, s \left(t \right),
\end{equation}
\begin{equation}\label{eq4.13}
H_5 \left(t \right)= \displaystyle \sum_{n=2}^{\infty} \left(h_n - g_n \right) \cos \left(n + \delta_n \left( \pi, \beta \right)\right)t.
\end{equation}
The absolute continuity of the series in \eqref{eq4.9} can be found in Lemma \ref{lem7.1}. The estimates for $s_n,$ $h_n,$ $g_n,$ $\epsilon_n,$ $m_n \left(t\right)$ and $z_n \left(t\right)$ ensure the absolutely and uniformly convergence of the series in \eqref{eq4.11}, \eqref{eq4.13} on $\left[0, 2\pi\right]$  and absolutely continuity of the functions $H_3$ and $H_5$ in $\left(0, 2\pi\right).$ \\
Since $l \left(t \right)$ and $s \left(t \right)$ are absolutely continuous functions in $\left(0, 2\pi \right),$ then Lemma \ref{lem4.2} is proved.
\end{proof}

\begin{remark}\label{rem4.3}
The function $H$ (see \eqref{eq4.3}) is defined for the case $\lambda_n \neq 0,$ $n=0,1,2,\dots$ and $\lambda_n \left(0, \pi, \beta \right) \neq 0,$ $n=0,1,2,\dots.$ For the other cases we define $H$ in the following way:
\begin{itemize}
  \item if there is $n=n_0,$ such that $\lambda_{n_0}=0$ and $\lambda_n \left(0, \pi, \beta \right) \neq 0,$ $n=0,1,2,\dots,$ then
\begin{multline*}
H \left(t\right)=-\dfrac{1}{a_{n_0}} \, \dfrac{t^2}{2}-\dfrac{1}{{a_{n_{0}}} \left( 0, \pi, \beta \right)} \, \dfrac{\cos \lambda_{n_{0}} \left(0, \pi, \beta \right) t} {\lambda_{n_{0}}^{2} \left(0, \pi, \beta \right)}+\\
+\sum\limits_{\substack{n = 0 \\ n \neq n_{0}}}^{\infty} \left(\dfrac{1}{a_n} \, \dfrac{\cos \lambda_{n} t}{\lambda_{n}^2}- \dfrac{1}{{a_n} \left( 0, \pi, \beta \right)} \, \dfrac{\cos \lambda_{n} \left(0, \pi, \beta \right) t} {\lambda_{n}^{2} \left(0, \pi, \beta \right)}\right),
\end{multline*}
  \item if $\lambda_n \neq 0,$ $n=0,1,2,\dots$ and there is $n=n_1,$ such that $\lambda_{n_1} \left(0, \pi, \beta \right)=0,$ then
\begin{multline*}
H \left(t\right)=\dfrac{1}{a_{n_1}} \, \dfrac{\cos \lambda_{n_1} t}{\lambda_{n_1}^2}+\dfrac{1}{{a_{n_1}} \left( 0, \pi, \beta \right)} \, \dfrac{t^2}{2}+\\
+\sum\limits_{\substack{n = 0 \\ n \neq n_{1}}}^{\infty} \left(\dfrac{1}{a_n} \, \dfrac{\cos \lambda_{n} t}{\lambda_{n}^2}- \dfrac{1}{{a_n} \left( 0, \pi, \beta \right)} \, \dfrac{\cos \lambda_{n} \left(0, \pi, \beta \right) t} {\lambda_{n}^{2} \left(0, \pi, \beta \right)}\right),
\end{multline*}
  \item if there are $n_0$ and $n_1,$ such that $\lambda_{n_0}=0$ and $\lambda_{n_1} \left(0, \pi, \beta \right)=0,$ then
\begin{multline*}
H \left(t\right)=\dfrac{1}{a_{n_1}} \, \dfrac{\cos \lambda_{n_1} t}{\lambda_{n_1}^2}-\dfrac{1}{{a_{n_{0}}} \left( 0, \pi, \beta \right)} \, \dfrac{\cos \lambda_{n_{0}} \left(0, \pi, \beta \right) t} {\lambda_{n_{0}}^{2} \left(0, \pi, \beta \right)}+\\
+\left(\dfrac{1}{{a_{n_{1}}} \left( 0, \pi, \beta \right)}-\dfrac{1}{a_{n_0}}\right) \dfrac{t^2}{2}+\sum\limits_{\substack{n = 0 \\ n \neq n_{0} \\ n \neq n_{1}}}^{\infty} \left(\dfrac{1}{a_n} \, \dfrac{\cos \lambda_{n} t}{\lambda_{n}^2}- \dfrac{1}{{a_n} \left( 0, \pi, \beta \right)} \, \dfrac{\cos \lambda_{n} \left(0, \pi, \beta \right) t} {\lambda_{n}^{2} \left(0, \pi, \beta \right)}\right).
\end{multline*}
\end{itemize}
It is easy to see that such defined functions are also absolutely continuous on arbitrary segment $\left[a, b \right] \subset \left(0, 2 \pi \right).$
\end{remark}
Consider the function
\begin{equation}\label{eq4.14}
F \left(x, t \right):=\dfrac{1}{2} \left( H \left( \left|x-t \right| \right)-H \left(x+t \right) \right).
\end{equation}
It follows from Lemma \ref{lem4.2} and Remark \ref{rem4.3} that,
\begin{itemize}
  \item if $\lambda_n \neq 0,$ $n=0,1,2,\dots$ and $\lambda_n \left(0, \pi, \beta \right) \neq 0,$ $n=0,1,2,\dots,$ then
\begin{equation}\label{eq4.15}
F\left(x, t \right)=\sum\limits_{n=0}^{\infty}\left(\dfrac{1}{a_n} \, \dfrac{\sin \lambda_{n} x}{\lambda_{n}} \, \dfrac{\sin \lambda_{n} t}{\lambda_{n}}- \dfrac{1}{{a_n} \left( 0, \pi, \beta \right)} \, \dfrac{\sin \lambda_{n} \left(0, \pi, \beta\right) x}{{\lambda_{n}\left(0, \pi, \beta \right)}} \, \dfrac{\sin \lambda_{n} \left(0, \pi, \beta \right) t}{{\lambda_{n}\left(0, \pi, \beta \right)}}\right),
\end{equation}
  \item if there is $n=n_0,$ such that $\lambda_{n_0}=0$ and $\lambda_n \left(0, \pi, \beta \right) \neq 0,$ $n=0,1,2,\dots,$ then
\begin{multline*}
F\left(x, t \right)=\dfrac{1}{a_{n_0}} \, x \, t-\dfrac{1}{{a_{n_0}} \left( 0, \pi, \beta \right)} \, \dfrac{\sin \lambda_{n_0} \left(0, \pi, \beta\right) x}{{\lambda_{n_0} \left(0, \pi, \beta \right)}} \, \dfrac{\sin \lambda_{n_0} \left(0, \pi, \beta \right) t}{{\lambda_{n_0} \left(0, \pi, \beta \right)}}+\\
+ \sum\limits_{\substack{n = 0 \\ n \neq n_{0}}}^{\infty} \left(\dfrac{1}{a_n} \, \dfrac{\sin \lambda_{n} x}{\lambda_{n}} \, \dfrac{\sin \lambda_{n} t}{\lambda_{n}}- \dfrac{1}{{a_n} \left( 0, \pi, \beta \right)} \, \dfrac{\sin \lambda_{n} \left(0, \pi, \beta\right) x}{{\lambda_{n}\left(0, \pi, \beta \right)}} \, \dfrac{\sin \lambda_{n} \left(0, \pi, \beta \right) t}{{\lambda_{n}\left(0, \pi, \beta \right)}}\right),
\end{multline*}
  \item if $\lambda_n \neq 0,$ $n=0,1,2,\dots$ and there is $n=n_1,$ such that $\lambda_{n_1} \left(0, \pi, \beta \right)=0,$ then
  \begin{multline*}
F\left(x, t \right)=\dfrac{1}{a_{n_1}} \, \dfrac{\sin \lambda_{n_1} x}{\lambda_{n_1}} \, \dfrac{\sin \lambda_{n_1} t}{\lambda_{n_1}}-\dfrac{1}{{a_{n_1}} \left( 0, \pi, \beta \right)} \, x \, t+\\
+ \sum\limits_{\substack{n = 0 \\ n \neq n_{1}}}^{\infty} \left(\dfrac{1}{a_n} \, \dfrac{\sin \lambda_{n} x}{\lambda_{n}} \, \dfrac{\sin \lambda_{n} t}{\lambda_{n}}- \dfrac{1}{{a_n} \left( 0, \pi, \beta \right)} \, \dfrac{\sin \lambda_{n} \left(0, \pi, \beta\right) x}{{\lambda_{n}\left(0, \pi, \beta \right)}} \, \dfrac{\sin \lambda_{n} \left(0, \pi, \beta \right) t}{{\lambda_{n}\left(0, \pi, \beta \right)}}\right),
\end{multline*}
  \item if there are $n_0$ and $n_1,$ such that $\lambda_{n_0}=0$ and $\lambda_{n_1} \left(0, \pi, \beta \right)=0,$ then
\begin{multline*}
F\left(x, t \right)=\dfrac{1}{a_{n_1}} \, \dfrac{\sin \lambda_{n_1} x}{\lambda_{n_1}} \, \dfrac{\sin \lambda_{n_1} t}{\lambda_{n_1}}+\left(\dfrac{1}{a_{n_0}}-\dfrac{1}{{a_{n_1}} \left( 0, \pi, \beta \right)}\right) x \, t-\\
-\dfrac{1}{{a_{n_0}} \left( 0, \pi, \beta \right)} \, \dfrac{\sin \lambda_{n_0} \left(0, \pi, \beta\right) x}{{\lambda_{n_0} \left(0, \pi, \beta \right)}} \, \dfrac{\sin \lambda_{n_0} \left(0, \pi, \beta \right) t}{{\lambda_{n_0} \left(0, \pi, \beta \right)}}+\\
+ \sum\limits_{\substack{n = 0 \\ n \neq n_{0} \\ n \neq n_{1}}}^{\infty} \left(\dfrac{1}{a_n} \, \dfrac{\sin \lambda_{n} x}{\lambda_{n}} \, \dfrac{\sin \lambda_{n} t}{\lambda_{n}}- \dfrac{1}{{a_n} \left( 0, \pi, \beta \right)} \, \dfrac{\sin \lambda_{n} \left(0, \pi, \beta\right) x}{{\lambda_{n}\left(0, \pi, \beta \right)}} \, \dfrac{\sin \lambda_{n} \left(0, \pi, \beta \right) t}{{\lambda_{n}\left(0, \pi, \beta \right)}}\right).
\end{multline*}
\end{itemize}
It is easy to see that $F\left(x, t \right)$ (in all four cases) is a continuous function and $\dfrac{d}{dx} F\left(x, x \right) \in L^1 \left(0, \pi \right).$

\begin{theorem}\label{thm4.4}
For each fixed $x \in \left(0, \pi \right],$ the kernel $P \left(x, t \right)$ of the transformation operator (see \eqref{eq4.1}) satisfies the following linear integral equation
\begin{equation}\label{eq4.16}
P \left(x, t \right)+F\left(x, t \right)+{\int_{0}^{x}}P\left(x, s \right) F\left(s, t \right)ds=0, \; 0 \leq t < x,
\end{equation}
which is also called Gelfand-Levitan equation.
\end{theorem}
\begin{proof}
We prove the theorem for the case $\lambda_n \neq 0,$ $n=0,1,2,\dots$ and $\lambda_n \left(0, \pi, \beta \right) \neq 0,$ $n=0,1,2,\dots,$ when the representation \eqref{eq4.15} holds. Other cases can be proven analogously.

Since $\mathbb{P}$ defined in \eqref{eq4.1} is a Volterra integral operator with continuous kernel $P(x,t),$ then $\mathbb{I}+\mathbb{P}$ has the inverse operator of the same type (see, for example, \cite{Mihlin:1949,Marchenko:1952,Marchenko:1977,Gasymov_Levitan:1964}), which we denote by $\mathbb{I}+\mathbb{Q}.$ Solving the equation \eqref{eq4.1} with respect to $\dfrac{\sin \lambda x}{\lambda}$ we obtain
\begin{equation}\label{eq4.17}
\dfrac{\sin \lambda x}{\lambda}=\left(\mathbb{I}+\mathbb{Q}\right)\varphi_{\pi} \left(x, {\lambda}^2 \right):=\varphi_{\pi} \left(x, {\lambda}^2 \right)+{\int_{0}^{x}}Q(x,t)\varphi_{\pi} \left(t, \lambda^2 \right)dt,
\end{equation}
where $Q\left(x, t \right),$ $0 \leq t \leq x \leq \pi$ is a real continuous function with the same smoothness as $P \left(x, t \right),$ and
\begin{equation}\label{eq4.18}
Q \left(x,x \right)=-\dfrac{1}{2}{\int_{0}^{x}}q\left(t \right)dt, \; Q \left(x, 0 \right)=0.
\end{equation}
In view of \eqref{eq4.1} and \eqref{eq4.17} we get that
\begin{multline}\label{eq4.19}
\sum\limits_{n = 2}^N \dfrac{\varphi_{\pi} \left(x, \lambda_n^2 \right){\dfrac{\sin {\lambda_n}t}{\lambda_n}}}{a_n}=\sum\limits_{n = 2}^N \dfrac{1}{a_n} \left( {\dfrac{\sin {\lambda_n}x}{\lambda_n}}+ {\int_{0}^{x}}P(x,s){\dfrac{\sin {\lambda_n}s}{\lambda_n}}ds \right){\dfrac{\sin {\lambda_n}t}{\lambda_n}}=\\
=\sum\limits_{n = 2}^N \dfrac{1}{a_n} \left( {\dfrac{\sin {\lambda_n}x \sin {\lambda_n}t}{\lambda_n^2}}+{\dfrac{\sin {\lambda_n}t}{\lambda_n}}{\int_{0}^{x}}P(x,s){\dfrac{\sin {\lambda_n}s}{\lambda_n}}ds \right)
\end{multline}
and
\begin{multline}\label{eq4.20}
\sum\limits_{n = 2}^N \dfrac{\varphi_{\pi} \left(x, \lambda_n^2 \right){\dfrac{\sin {\lambda_n}t}{\lambda_n}}}{a_n}=\sum\limits_{n = 2}^N \dfrac{1}{a_n} \varphi_{\pi} \left(x, \lambda_n^2 \right) \left( \varphi_{\pi} \left(t, \lambda_n^2 \right)+ {\int_{0}^{t}}Q(t,s)\varphi_{\pi} \left(s, \lambda_n^2 \right)ds \right)=\\
=\sum\limits_{n = 2}^N \dfrac{1}{a_n} \left( \varphi_{\pi} \left(x, \lambda_n^2 \right) \varphi_{\pi} \left(t, \lambda_n^2 \right) + \varphi_{\pi} \left(x, \lambda_n^2 \right){\int_{0}^{t}}Q(t,s)\varphi_{\pi} \left(s, \lambda_n^2 \right)ds \right).
\end{multline}
By $\Phi_N \left(x, t \right)$ denote the following expression
\begin{multline}\label{eq4.21}
\Phi_N \left(x, t \right):=\sum\limits_{n=0}^1 \left( \dfrac{1}{a_n} \, \varphi_{\pi} \left(x, \lambda_n^2 \right)\varphi_{\pi} \left(t, \lambda_n^2 \right)-\dfrac{1}{a_n \left(0, \pi, \beta \right)} \, \dfrac{\sin \lambda_n \left(0, \pi, \beta\right)x}{\lambda_n \left(0, \pi, \beta\right)} \dfrac{\sin \lambda_n \left(0, \pi, \beta\right)t}{\lambda_n \left(0, \pi, \beta\right)} \right)+\\
+\sum\limits_{n=2}^N \left( \dfrac{1}{a_n} \, \varphi_{\pi} \left(x, \lambda_n^2 \right)\varphi_{\pi} \left(t, \lambda_n^2 \right)-\dfrac{1}{a_n \left(0, \pi, \beta \right)} \, \dfrac{\sin \left(n+\delta_n \left(\pi, \beta \right)\right)x}{n+\delta_n \left(\pi, \beta \right)} \, \dfrac{\sin \left(n+\delta_n \left(\pi, \beta \right)\right)t}{n+\delta_n \left(\pi, \beta \right)} \right).
\end{multline}
It is easy to calculate that the following representation holds (compare with \cite{Freiling_Yurko:2001})
\begin{equation}\label{eq4.22}
\Phi_N \left(x, t \right)=I_{N1} \left(x, t \right)+I_{N2} \left(x, t \right)+I_{N3} \left(x, t \right)+I_{N4} \left(x, t \right),
\end{equation}
where
\begin{multline}\label{eq4.23}
I_{N1} \left(x, t \right)=\sum\limits_{n=0}^1 \left(\dfrac{1}{a_n} \, \dfrac{\sin {\lambda_n}x}{\lambda_n} \, \dfrac{\sin {\lambda_n}t}{\lambda_n}-\dfrac{1}{a_n \left(0, \pi, \beta \right)} \, \dfrac{\sin \lambda_n \left(0, \pi, \beta\right)x}{\lambda_n \left(0, \pi, \beta\right)} \, \dfrac{\sin \lambda_n \left(0, \pi, \beta\right)t}{\lambda_n \left(0, \pi, \beta\right)} \right)+\\
+\sum\limits_{n=2}^N \left(\dfrac{1}{a_n} \, \dfrac{\sin {\lambda_n}x}{\lambda_n} \, \dfrac{\sin {\lambda_n}t}{\lambda_n}-\dfrac{1}{a_n \left(0, \pi, \beta \right)} \, \dfrac{\sin \left(n+\delta_n(\pi, \beta)\right)x}{n+\delta_n(\pi, \beta)} \, \dfrac{\sin \left(n+\delta_n(\pi, \beta)\right)t}{n+\delta_n(\pi, \beta)} \right),
\end{multline}
\begin{multline}\label{eq4.24}
I_{N2} \left(x, t \right)=\sum\limits_{n=0}^1 \dfrac{1}{a_n \left(0, \pi, \beta \right)} \, \dfrac{\sin \lambda_n \left(0, \pi, \beta\right)t}{\lambda_n \left(0, \pi, \beta\right)}{\int_{0}^{x}}P(x,s)\dfrac{\sin \lambda_n \left(0, \pi, \beta\right)s}{\lambda_n \left(0, \pi, \beta\right)}ds+\\
+\sum\limits_{n=2}^N \dfrac{1}{a_n \left(0, \pi, \beta \right)} \, \dfrac{\sin \left(n+\delta_n(\pi, \beta)\right)t}{n+\delta_n(\pi, \beta)}{\int_{0}^{x}}P(x,s)\dfrac{\sin \left(n+\delta_n(\pi, \beta)\right)s}{n+\delta_n(\pi, \beta)}ds,
\end{multline}
\begin{multline}\label{eq4.25}
I_{N3} \left(x, t \right)=\sum\limits_{n=0}^1 {\int_{0}^{x}}P \left(x, s \right) \left(\dfrac{1}{a_n} \, \dfrac{\sin {\lambda_n}s}{\lambda_n} \, \dfrac{\sin {\lambda_n}t}{\lambda_n}\right.-\\
-\left.\dfrac{1}{a_n \left(0, \pi, \beta \right)} \, \dfrac{\sin \lambda_n \left(0, \pi, \beta\right)s}{\lambda_n \left(0, \pi, \beta\right)} \, \dfrac{\sin \lambda_n \left(0, \pi, \beta\right)t}{\lambda_n \left(0, \pi, \beta\right)} \right)ds+\\
+\sum\limits_{n=2}^N {\int_{0}^{x}}P(x,s) \left(\dfrac{1}{a_n} \, \dfrac{\sin {\lambda_n}s}{\lambda_n} \, \dfrac{\sin {\lambda_n}t}{\lambda_n}\right.-\\
-\left.\dfrac{1}{a_n \left(0, \pi, \beta \right)} \, \dfrac{\sin \left(n+\delta_n(\pi, \beta)\right)s}{n+\delta_n(\pi, \beta)} \, \dfrac{\sin \left(n+\delta_n(\pi, \beta)\right)t}{n+\delta_n(\pi, \beta)} \right)ds,
\end{multline}
\begin{equation}\label{eq4.26}
I_{N4} \left(x, t \right)=-\sum\limits_{n=0}^N \dfrac{\varphi_{\pi} \left(x, \lambda_n^2 \right)}{a_n}{\int_{0}^{t}}Q(t,s) \varphi_{\pi} \left(s, \lambda_n^2 \right)ds.
\end{equation}
Let $f$ is an absolutely continuous function on $\left[0, \pi \right],$ such that $f\left(0\right)=0.$ In light of Theorem \ref{thm2.3} and Remark \ref{rem2.7} we receive
\begin{multline}\label{eq4.27}
\lim\limits_{N \to \infty} \mathop {\operatorname{max}} \limits_{x \in \left[0, \pi\right]} \displaystyle\int_{0}^{\pi}f\left(t\right) \Phi_N\left(x,t\right)dt =\lim\limits_{N \to \infty} \mathop {\operatorname{max}} \limits_{x \in \left[0, \pi\right]} \sum\limits_{n=0}^N \dfrac{1}{a_n}\int_{0}^{\pi}f\left(t\right)\varphi_{\pi} \left(t, \lambda_n^2 \right)dt \varphi_{\pi} \left(x, \lambda_n^2 \right)-\\
-\lim\limits_{N \to \infty} \mathop {\operatorname{max}} \limits_{x \in \left[0, \pi\right]} \sum\limits_{n=0}^N \dfrac{1}{a_n \left(0, \pi, \beta \right)}\int_{0}^{\pi}f\left(t\right)\dfrac{\sin \lambda_n \left(0, \pi, \beta\right)t}{\lambda_n \left(0, \pi, \beta\right)}dt \dfrac{\sin \lambda_n \left(0, \pi, \beta\right)x}{\lambda_n \left(0, \pi, \beta\right)}=f\left(x\right)-f\left(x\right)=0.
\end{multline}
Extending $P\left(x, t \right)=Q\left(x, t \right)=0$ for $x<t$ and using the standard theorems of the limit process under the integral sign we get that uniformly with respect to $x \in \left[0, \pi \right]$
\begin{multline}\label{eq4.28}
\lim\limits_{N \to \infty} \displaystyle\int_{0}^{\pi}f\left(t\right) I_{N1}\left(x,t\right)dt =\lim\limits_{N \to \infty} \displaystyle\int_{0}^{\pi}f\left(t\right) \sum\limits_{n=0}^N \left(\dfrac{1}{a_n} \dfrac{\sin {\lambda_n}x}{\lambda_n} \dfrac{\sin {\lambda_n}t}{\lambda_n}\right.-\\
-\left.\dfrac{1}{a_n \left(0, \pi, \beta \right)}\dfrac{\sin \lambda_n \left(0, \pi, \beta\right)x}{\lambda_n \left(0, \pi, \beta\right)} \dfrac{\sin \lambda_n \left(0, \pi, \beta\right)t}{\lambda_n \left(0, \pi, \beta\right)} \right)dt=\displaystyle\int_{0}^{\pi}f\left(t\right) \lim\limits_{N \to \infty} \sum\limits_{n=0}^N \left(\dfrac{1}{a_n} \dfrac{\sin {\lambda_n}x}{\lambda_n} \dfrac{\sin {\lambda_n}t}{\lambda_n}\right.-\\
-\left.\dfrac{1}{a_n \left(0, \pi, \beta \right)}\dfrac{\sin \lambda_n \left(0, \pi, \beta\right)x}{\lambda_n \left(0, \pi, \beta\right)} \dfrac{\sin \lambda_n \left(0, \pi, \beta\right)t}{\lambda_n \left(0, \pi, \beta\right)} \right)dt=\displaystyle\int_{0}^{\pi}f\left(t\right)F\left(x,t\right)dt,
\end{multline}
\begin{multline}\label{eq4.29}
\lim\limits_{N \to \infty} \displaystyle\int_{0}^{\pi}f\left(t\right) I_{N2}\left(x,t\right)dt =\\
=\lim\limits_{N \to \infty} \displaystyle\int_{0}^{\pi}f\left(t\right) \sum\limits_{n=0}^N \dfrac{1}{a_n \left(0, \pi, \beta \right)} \dfrac{\sin \lambda_n \left(0, \pi, \beta\right)t}{\lambda_n \left(0, \pi, \beta\right)}{\int_{0}^{x}}P(x,s)\dfrac{\sin \lambda_n \left(0, \pi, \beta\right)s}{\lambda_n \left(0, \pi, \beta\right)}dsdt=\\
=\lim\limits_{N \to \infty} \sum\limits_{n=0}^N \dfrac{1}{a_n \left(0, \pi, \beta \right)} \displaystyle\int_{0}^{\pi} \displaystyle\int_{0}^{\pi} f\left(t\right)  \dfrac{\sin \lambda_n \left(0, \pi, \beta\right)t}{\lambda_n \left(0, \pi, \beta\right)} P(x,s) \dfrac{\sin \lambda_n \left(0, \pi, \beta\right)s}{\lambda_n \left(0, \pi, \beta\right)} dtds=\\
=\displaystyle\int_{0}^{\pi} \lim\limits_{N \to \infty} \sum\limits_{n=0}^N \dfrac{1}{a_n \left(0, \pi, \beta \right)} \displaystyle\int_{0}^{\pi}f\left(t\right)  \dfrac{\sin \lambda_n \left(0, \pi, \beta\right)t}{\lambda_n \left(0, \pi, \beta\right)}dt\dfrac{\sin \lambda_n \left(0, \pi, \beta\right)s}{\lambda_n \left(0, \pi, \beta\right)} P(x,s)ds =\\
=\displaystyle\int_{0}^{\pi}f\left(s\right) P(x,s)ds=\displaystyle\int_{0}^{\pi}f\left(t\right) P(x,t)dt,
\end{multline}
\begin{multline}\label{eq4.30}
\lim\limits_{N \to \infty} \displaystyle\int_{0}^{\pi}f\left(t\right) I_{N3}\left(x,t\right)dt=\lim\limits_{N \to \infty} \displaystyle\int_{0}^{\pi}f\left(t\right) \sum\limits_{n=0}^N {\int_{0}^{x}}P(x,s) \left(\dfrac{1}{a_n} \dfrac{\sin {\lambda_n}s}{\lambda_n} \dfrac{\sin {\lambda_n}t}{\lambda_n}\right.-\\
-\left.\dfrac{1}{a_n \left(0, \pi, \beta \right)}\dfrac{\sin \lambda_n \left(0, \pi, \beta\right)s}{\lambda_n \left(0, \pi, \beta\right)} \dfrac{\sin \lambda_n \left(0, \pi, \beta\right)t}{\lambda_n \left(0, \pi, \beta\right)} \right)dsdt=\displaystyle\int_{0}^{\pi}f\left(t\right) {\int_{0}^{x}}P(x,s) \times\\
\times \lim\limits_{N \to \infty}\sum\limits_{n=0}^N \left(\dfrac{1}{a_n} \dfrac{\sin {\lambda_n}s}{\lambda_n} \dfrac{\sin {\lambda_n}t}{\lambda_n}-\dfrac{1}{a_n \left(0, \pi, \beta \right)}\dfrac{\sin \lambda_n \left(0, \pi, \beta\right)s}{\lambda_n \left(0, \pi, \beta\right)} \dfrac{\sin \lambda_n \left(0, \pi, \beta\right)t}{\lambda_n \left(0, \pi, \beta\right)} \right)dsdt=\\
=\displaystyle\int_{0}^{\pi}f\left(t\right)\left( {\int_{0}^{x}}P(x,s)F(s,t)ds\right)dt,
\end{multline}
\begin{multline}\label{eq4.31}
\lim\limits_{N \to \infty} \displaystyle\int_{0}^{\pi}f\left(t\right) I_{N4}\left(x,t\right)dt=-\lim\limits_{N \to \infty} \displaystyle\int_{0}^{\pi}f\left(t\right)\sum\limits_{n=0}^N \dfrac{\varphi_{\pi} \left(x, \lambda_n^2 \right)}{a_n}{\int_{0}^{t}}Q(t,s) \varphi_{\pi} \left(s, \lambda_n^2 \right)dsdt=\\
=-\lim\limits_{N \to \infty} \sum\limits_{n=0}^N \dfrac{\varphi_{\pi} \left(x, \lambda_n^2 \right)}{a_n}\displaystyle\int_{0}^{\pi}f\left(t\right){\int_{0}^{t}}Q(t,s) \varphi_{\pi} \left(s, \lambda_n^2 \right)dsdt=\\
=-\lim\limits_{N \to \infty} \sum\limits_{n=0}^N \dfrac{\varphi_{\pi} \left(x, \lambda_n^2 \right)}{a_n}\displaystyle\int_{0}^{\pi}\left({\int_{s}^{\pi}}Q(t,s) f\left(t\right)dt\right)\varphi_{\pi} \left(s, \lambda_n^2 \right)ds=-{\int_{0}^{\pi}}Q(t,x) f\left(t\right)dt ,
\end{multline}
where $Q\left(t, x \right)=0$ for $t<x.$
From \eqref{eq4.22} and \eqref{eq4.27}--\eqref{eq4.31}, we can conclude that
\begin{equation}\label{eq4.32}
\displaystyle\int_{0}^{\pi} f\left(t\right) \left(F\left(x,t\right)+P\left(x,t\right)+{\int_{0}^{x}}P(x,s)F(s,t)ds-Q\left(t,x\right)\right)dt=0.
\end{equation}
Since the system of eigenfunctions $\left\{\varphi_{m} \left(t \right)\right\}_{m=0}^\infty$ of the boundary value problem $L(q, \pi, \beta)$ is complete in $L^2 \left(0, \pi \right)$ and $\varphi_{m} \left(0 \right)=0,$ $m=0,1,2,\dots,$ then we can take $f\left(t \right)=\varphi_{m} \left(t \right),$ $m=0,1,2,\dots$ and obtain that for each fixed $x \in \left(0, \pi \right]$
\begin{equation}\label{eq4.33}
F\left(x,t\right)+P\left(x,t\right)+\displaystyle\int_{0}^{x} P(x,s)F(s,t)ds-Q\left(t,x\right)=0.
\end{equation}
For $t<x,$ this is equivalent to \eqref{eq4.16}. This completes the proof.
\end{proof}

\begin{remark}\label{rem4.5}
It follows from \eqref{eq4.2}, \eqref{eq4.18} and \eqref{eq4.33} that for all $x \in \left[0, \pi \right]$
\begin{equation}\label{eq4.34}
2 P\left(x, x\right)+F\left(x, x\right)+\displaystyle\int_{0}^{x} P(x,s)F(s,x)ds=0.
\end{equation}
\end{remark}

\section{The constructive solution of the inverse problem}
\label{sec5}

\begin{lemma}\label{lem5.1}
For each fixed $x \in \left( 0, \pi \right],$ equation \eqref{eq4.16} has a unique solution $P\left(x, \cdot \right)$ in $L^2\left[0, x \right).$
\end{lemma}
\begin{proof}
Since \eqref{eq4.16} is a Fredholm equation, then it is sufficient (see, for example, \cite{Mihlin:1949}) to prove that the homogenous equation
\begin{equation}\label{eq5.1}
p \left(t \right)+{\int_{0}^{x}}p \left(s \right)F \left(s, t \right)ds=0
\end{equation}
has only a trivial solution $p \left(t \right)=0.$

Without loss of generality (see Remark \ref{rem5.2}) we suppose that $\lambda_n \neq 0,$ $n=0,1,2,\dots$ and $\lambda_n \left(0, \pi, \beta \right) \neq 0,$ $n=0,1,2,\dots.$ Let $p(t)$ be a solution of \eqref{eq5.1}. Then multiplying \eqref{eq5.1} by $p\left(t\right)$ and integrating from $0$ to $x,$ we receive
\begin{equation*}
{\int_{0}^{x}}{p^2} \left(t \right)dt+{\int_{0}^{x}}{\int_{0}^{x}}F \left(s, t \right)p \left(s \right)p \left(t \right)dsdt=0
\end{equation*}
or
\begin{multline*}
{\int_{0}^{x}}{p^2}(t)dt+\sum\limits_{n = 0}^\infty \dfrac{1}{a_n} {\int_{0}^{x}} p \left(s \right) \dfrac{\sin {\lambda_n}s}{\lambda_n}ds{\int_{0}^{x}} p \left(t \right) \dfrac{\sin {\lambda_n}t}{\lambda_n}dt-\\
-\sum\limits_{n = 0}^\infty \dfrac{1}{a_n \left(0, \pi, \beta \right)} {\int_{0}^{x}} p \left(s \right) \dfrac{\sin {\lambda_n \left(0, \pi, \beta\right)}s}{\lambda_n \left(0, \pi, \beta\right)}ds {\int_{0}^{x}} p \left(t \right) \dfrac{\sin {\lambda_n \left(0, \pi, \beta\right)}t}{\lambda_n \left(0, \pi, \beta\right)}dt=0,
\end{multline*}
which is equivalent to
\begin{multline*}
{\int_{0}^{x}}{p^2} \left(t \right)dt+\sum\limits_{n = 0}^\infty \dfrac{1}{a_n} \left({\int_{0}^{x}} p \left(s \right) \dfrac{\sin {\lambda_n}t}{\lambda_n}dt\right)^2-\\
-\sum\limits_{n = 0}^\infty \dfrac{1}{a_n \left(0, \pi, \beta \right)} \left( {\int_{0}^{x}} p \left(t \right) \dfrac{\sin \lambda_n \left(0, \pi, \beta\right)t}{\lambda_n \left(0, \pi, \beta\right)}dt\right)^2=0.
\end{multline*}
Since $\left\{\dfrac{\sin \lambda_n \left(0, \pi, \beta\right)t}{\lambda_n \left(0, \pi, \beta\right)}\right\}_{n=0}^\infty$ are the eigenfunctions of the problem $L \left(0, \pi, \beta \right),$ then we can extend the function $p \left(t \right)$ by zero for $t>x$ and write Parseval's identity in the following form
\begin{equation*}
{\int_{0}^{\pi}}{p^2} \left(t \right)dt={\int_{0}^{x}}{p^2} \left(t \right)dt=\sum\limits_{n = 0}^\infty \dfrac{1}{a_n \left(0, \pi, \beta \right)} \left( {\int_{0}^{x}} p \left(t \right) \dfrac{\sin \lambda_n \left(0, \pi, \beta\right)t}{\lambda_n \left(0, \pi, \beta\right)}dt\right)^2,
\end{equation*}
and therefore,
\begin{equation*}
\sum\limits_{n = 0}^\infty \dfrac{1}{a_n} \left({\int_{0}^{x}} p \left(s \right) \dfrac{\sin {\lambda_n}t}{\lambda_n}dt\right)^2=0.
\end{equation*}
Since $a_n>0$ (see \eqref{eq1.12}), then
\begin{equation*}
{\int_{0}^{x}} p \left(s \right) \sin {\lambda_n}tdt=0, \; n \geq 0.
\end{equation*}
In \cite[Theorem 1.2, proposition 2a]{Harutyunyan_Pahlevanyan_Srapionyan:2013} it was proved that the system of functions $\left\{\sin {\lambda_n}t\right\}_{n=0}^\infty,$ where $\lambda_n=\lambda_n \left(q, \pi, \beta \right) \neq 0,$ is a Riesz basis in  $L^2\left[0, \pi\right].$ In fact, there was proved a stronger statement: if the sequence $\lambda_n \neq 0,$ $n=0,1,2,\dots$ has the following asymptotic property:
\begin{equation}\label{eq5.2}
\lambda_n=n+\delta_n \left(\pi,\beta\right)+O \left(n^{-1}\right),
\end{equation}
then $\left\{\sin {\lambda_n}t\right\}_{n=0}^\infty$ is a Riesz basis in  $L^2\left[0, \pi\right].$

Since our given numbers $\left\{\lambda_n\right\}_{n \geq 0}$ (see \eqref{eq1.11}) satisfy asymptotic property \eqref{eq5.2}, then $\left\{\sin {\lambda_n}t\right\}_{n=0}^\infty$ is a Riesz basis and accordingly is complete in $L^2\left[0, \pi\right].$ Hence, we have $p\left(t\right)=0.$ This completes the proof.
\end{proof}

Observe that the solution $P \left(x, t \right)$ of \eqref{eq4.16} has the same smoothness as $F \left(x, t \right)$ (see, e.g. \cite[equation (1.5.16) on p. 39]{Freiling_Yurko:2001}).

\begin{remark}\label{rem5.2}
If there are $n_1$ and $n_0$ such that $\lambda_{n_1}=0$ and (or) $\lambda_{n_0} \left(0, \pi, \beta \right)=0,$ then in the proof of the Lemma \ref{lem5.1} instead of $\left\{\dfrac{\sin \lambda_n \left(0, \pi, \beta\right)t}{\lambda_n \left(0, \pi, \beta\right)}\right\}_{n=0}^\infty$ we should take $\left\{\dfrac{\sin \lambda_n \left(0, \pi, \beta\right)t}{\lambda_n \left(0, \pi, \beta\right)}\right\}_{n=0}^{{n_0}-1} \cup \{t\} \cup \left\{\dfrac{\sin \lambda_n \left(0, \pi, \beta\right)t}{\lambda_n \left(0, \pi, \beta\right)}\right\}_{{n_0}+1}^\infty$ eigenfunctions of the problem $L \left(0, \pi, \beta \right)$ and instead of $\left\{\sin {\lambda_n}t\right\}_{n=0}^\infty$ we should take $\left\{\sin {\lambda_n}t\right\}_{n=0}^{{n_1}-1} \cup \{t\} \cup \left\{\sin {\lambda_n}t\right\}_{{n_1}+1}^\infty,$ which is a Riesz basis in $L^2\left[0, \pi\right]$ (see \cite[Theorem 1.2, proposition 2b]{Harutyunyan_Pahlevanyan_Srapionyan:2013}).
\end{remark}

\begin{remark}\label{rem5.3}
It is noteworthy that, in contrast to the above-mentioned classical results (see \cite{Gelfand_Levitan:1951, Gasymov_Levitan:1964,Zhikov:1967,Freiling_Yurko:2001}), in our case the Levinson's theorem (see \cite{Levinson:1940,Young:1980,Kadec:1964}) cannot be applied for proving completeness of the system of functions $\left\{\sin {\lambda_n}t\right\}_{n=0}^\infty,$ because of asymptotic behavior of $\lambda_n.$
\end{remark}
Let us define
\begin{equation}\label{eq5.3}
\varphi_{\pi}\left(x, \mu\right)\equiv \varphi_{\pi}\left(x, \lambda^2\right):=\dfrac{\sin \lambda x}{\lambda}+{\displaystyle\int_{0}^{x}}P \left(x, t \right)\dfrac{\sin \lambda t} {\lambda}dt
\end{equation}
and
\begin{equation}\label{eq5.4}
q\left(x \right):=2\dfrac{d}{dx}P\left(x, x \right).
\end{equation}
\begin{lemma}\label{lem5.4}
The following relations hold
\begin{equation}\label{eq5.5}
-\varphi_{\pi}''\left(x, \mu \right)+q\left(x \right)\varphi_{\pi} \left(x, \mu \right)=\mu \varphi_{\pi} \left(x, \mu \right),
\end{equation}
\begin{equation}\label{eq5.6}
\varphi_{\pi}\left(0, \mu \right)=0, \; \varphi_{\pi}'\left(0, \mu \right)=1.
\end{equation}
\end{lemma}
\begin{proof}
We split the proof into two cases:

{\bf Case I:} Let $H'\left(\cdot\right) \in AC\left(0, 2\pi\right),$ where $H\left(t\right)$ defined by \eqref{eq4.3} (here, again, we prove the lemma only for the case $\lambda_n \neq 0,$ $n=0,1,2,\dots$ and $\lambda_n \left(0, \pi, \beta \right) \neq 0,$ $n=0,1,2,\dots,$ other cases can be proven analogously). Differentiating the identity
\begin{equation}\label{eq5.7}
J\left(x, t \right) \equiv P\left(x, t \right)+F\left(x, t \right)+{\int_{0}^{x}}P\left(x, s \right) F\left(s, t \right) ds \equiv 0,
\end{equation}
we get
\begin{equation}\label{eq5.8}
J_{t}\left(x, t \right)=P_{t}\left(x, t \right)+F_{t}\left(x, t \right)+{\int_{0}^{x}}P\left(x, s \right) F_{t}\left(s, t\right)ds=0,
\end{equation}
\begin{equation}\label{eq5.9}
J_{tt}\left(x ,t \right)=P_{tt}\left(x, t \right)+F_{tt}\left(x, t \right)+{\int_{0}^{x}}P\left(x, s \right) F_{tt}\left(s, t \right)ds=0,
\end{equation}
\begin{equation}\label{eq5.10}
J_{x}\left(x, t \right)=P_{x}\left(x, t \right)+F_{x}\left(x, t \right)+P\left(x, x \right)F\left(x, t\right)+ {\int_{0}^{x}}P_{x}\left(x, s\right) F\left(s, t\right)ds=0,
\end{equation}
\begin{multline}\label{eq5.11}
J_{xx}\left(x, t\right)=P_{xx}\left(x, t\right)+F_{xx}\left(x, t\right)+\dfrac{dP\left(x, x\right)}{dx} F\left(x, t\right)+P\left(x, x\right)F_{x}\left(x, t\right)+\\
+\left.\dfrac{\partial P\left(x, s\right)}{\partial x}\right|_{s=x} F\left(x, t\right)+{\int_{0}^{x}}P_{xx}\left(x, s \right) F\left(s, t \right)ds=0.
\end{multline}
It is easy to see that $F_{xx}\left(x, t \right)=F_{tt}\left(x, t \right),$ $F\left(0, t \right)=0,$ $F\left(x, 0 \right)=0.$ Then, the identity \eqref{eq5.7} for $t=0$ gives $P\left(x, 0 \right)=0.$

After integrating \eqref{eq5.9} by parts twice, we obtain
\begin{multline}\label{eq5.12}
J_{tt}\left(x, t \right)=P_{tt}\left(x, t \right)+F_{tt}\left(x, t \right)+P\left(x, x \right) \left.\dfrac{\partial F\left(s, t \right)}{\partial s}\right|_{s=x}-\\
-F\left(x, t \right)\left.\dfrac{\partial P\left(x, s\right)}{\partial s} \right|_{s=x}+{\int_{0}^{x}}P_{ss}\left(x, s \right) F\left(s, t \right)ds=0.
\end{multline}
Consider the following identity
\begin{multline}\label{eq5.13}
J_{xx}\left(x, t \right)-J_{tt}\left(x, t \right)-q\left(x \right)J\left(x, t \right) \equiv P_{xx}\left(x, t \right)+\dfrac{dP\left(x, x \right)}{dx} F\left(x, t \right)+P\left(x, x \right)F_{x}\left(x, t\right)+\\
+\left.\dfrac{\partial P\left(x, s\right)}{\partial x}\right|_{s=x}F\left(x, t \right) +{\int_{0}^{x}}P_{xx}\left(x, s \right) F\left(s, t \right)ds-P_{tt}\left(x, t\right)-\\
-P\left(x, x\right)\left.\dfrac{\partial F\left(s, t\right)}{\partial s}\right|_{s=x}+F\left(x, t\right)\left.\dfrac{\partial P\left(x, s\right)}{\partial s}\right|_{s=x}-{\int_{0}^{x}}P_{ss}\left(x, s \right) F\left(s, t \right)ds-\\
-q\left(x \right)\left(P\left(x, t \right)+F\left(x, t \right)+{\int_{0}^{x}}P\left(x, s \right)F\left(s, t \right)ds \right) \equiv 0.
\end{multline}
Since $q\left(x\right) \equiv 2\dfrac{dP\left(x,x\right)}{dx} \equiv 2\left(\left.\dfrac{\partial P\left(x, s \right)}{\partial x}\right|_{s=x}+\left.\dfrac{\partial P\left(x, s \right)}{\partial s}\right|_{s=x}\right),$ then \eqref{eq5.13} is equivalent to the identity
\begin{multline}\label{eq5.14}
J_{xx}\left(x, t \right)-J_{tt}\left(x, t \right)-q\left(x \right)J\left(x, t \right) \equiv P_{xx}\left(x, t \right)-P_{tt}\left(x, t \right)-q\left(x \right)P\left(x, t \right)+\\
+{\int_{0}^{x}}\left(P_{xx}\left(x, s \right)-P_{ss}\left(x, s \right)-q\left(x \right)P\left(x, s \right)\right) F\left(s, t \right)ds \equiv 0.
\end{multline}
According to Lemma \ref{lem5.1}, the last equation has only the trivial solution, i.e.
\begin{equation}\label{eq5.15}
P_{xx}\left(x, t \right)-P_{tt}\left(x, t \right)-q\left(x \right)P\left(x, t \right)=0, \; 0 \leq t < x.
\end{equation}
On one hand, differentiating $\varphi_{\pi}\left(x, \mu \right)$ (see \eqref{eq5.3}) twice, we obtain
\begin{equation}\label{eq5.16}
\varphi_{\pi}'\left(x, \mu \right)=\cos \lambda x+P\left(x, x \right)\dfrac{\sin \lambda x}{\lambda}+{\int_{0}^{x}}P_{x}\left(x, t \right)\dfrac{\sin \lambda t}{\lambda}dt
\end{equation}
and
\begin{multline}\label{eq5.17}
\varphi_{\pi}''\left(x, \mu \right)=-\lambda \sin \lambda x+\dfrac{dP\left(x, x \right)}{dx}\dfrac{\sin \lambda x}{\lambda}+P\left(x, x \right) \cos \lambda x+\\
+\left.\dfrac{\partial P\left(x, t \right)}{\partial x}\right|_{t=x}\dfrac{\sin \lambda x}{\lambda}+{\int_{0}^{x}}P_{xx}\left(x, t \right)\dfrac{\sin \lambda t}{\lambda}dt.
\end{multline}
On the other hand, integrating $\varphi_{\pi}\left(x,\mu\right)$ by parts twice, we get
\begin{equation}\label{eq5.18}
\mu \varphi_{\pi}\left(x, \mu\right)=\lambda \sin \lambda x - P\left(x, x \right)\cos \lambda x + \left.\dfrac{\partial P\left(x, t \right)}{\partial t}\right|_{t=x}\dfrac{\sin \lambda x}{\lambda}- {\int_{0}^{x}}P_{tt}\left(x, t \right) \dfrac{\sin \lambda t}{\lambda}dt.
\end{equation}
Together with \eqref{eq5.17} and \eqref{eq5.15} this ensures that \eqref{eq5.5} holds and $\varphi_{\pi}\left(0, \mu \right)=0.$ From \eqref{eq5.16} it follows, that $\varphi_{\pi}'\left(0, \mu \right)=1,$ therefore \eqref{eq5.6} also holds.

{\bf Case II:} Now, let us consider the general case, when \eqref{eq1.11}--\eqref{eq1.14} hold and, according to Lemma \ref{lem4.2}, $H\left(\cdot\right) \in AC\left(0, 2\pi\right).$ Denote by $\tilde \varphi \left(x, \mu \right)$ solution of the equation \eqref{eq1.1} (where $q \left(x \right)=2\dfrac{d}{dx}P \left(x, x \right)$) with initial conditions $\tilde \varphi \left(0, \mu \right)=0,$  $\tilde \varphi' \left(0, \mu \right)=1.$ Our goal is to prove that $\tilde \varphi \left(x, \mu \right) \equiv \varphi_{\pi} \left(x, \mu \right).$

Choose the numbers $\left\{\lambda_{n,\left(j\right)}\right\}_{n \geq 0},$ $\left\{a_{n,\left(j\right)}\right\}_{n \geq 0},$ $j \geq 1$ of the form
\begin{equation}\label{eq5.19}
\lambda_{n,\left(j\right)}=n+\delta_n \left(\pi, \beta \right)+\dfrac {c}{2 \left(n+\delta_n \left(\pi, \beta \right) \right)}+l_{n,\left(j\right)}, \; \mu_{n,\left(j\right)} \neq \mu_{m,\left(j\right)} \, \left(n \neq m \right),
\end{equation}
\begin{equation}\label{eq5.20}
a_{n,\left(j\right)}=\dfrac {\pi}{2 \left(n+{\delta_n} \left(\pi, \beta \right) \right)^2} \left(1+\dfrac{c_1}{2 \left(n+\delta_n \left(\pi, \beta \right) \right)^2} +\dfrac{2 \, s_{n,\left(j\right)}}{\pi \left(n+\delta_n \left(\pi, \beta \right) \right)^2}\right), \; a_{n,\left(j\right)}>0,
\end{equation}
where $c_1$ is a constant, $l_{n,\left(j\right)}=o\left(\dfrac{1}{n^2}\right)$ and $s_{n,\left(j\right)}=o\left(1\right),$ are such that the functions
\begin{equation}\label{eq5.21}
l_j \left(t\right)=\displaystyle\sum_{n=2}^{\infty} l_{n,\left(j\right)} \sin \left(n+\delta_n \left(\pi, \beta\right)\right)t,
\end{equation}
\begin{equation}\label{eq5.22}
l'_j \left(t\right)=\displaystyle\sum_{n=2}^{\infty} l_{n,\left(j\right)} \left(n+\delta_n \left(\pi, \beta\right)\right) \sin \left(n+\delta_n \left(\pi, \beta\right)\right)t,
\end{equation}
\begin{equation}\label{eq5.23}
s_j \left(t\right)=\displaystyle \sum_{n=2}^{\infty} \dfrac {s_{n,\left(j\right)}}{\left(n + \delta_n \left( \pi, \beta \right)\right)^2} \cos \left(n + \delta_n \left( \pi, \beta \right)\right)t,
\end{equation}
\begin{equation}\label{eq5.24}
s'_j \left(t\right)=-\displaystyle \sum_{n=2}^{\infty} \dfrac {s_{n,\left(j\right)}}{n + \delta_n \left( \pi, \beta \right)} \sin \left(n + \delta_n \left( \pi, \beta \right)\right)t
\end{equation}
are absolutely continuous on arbitrary segment $\left[a, b \right] \subset (0, 2 \pi)$ (compare with asymtotics of eigenvalues and norming constants from Theorem \ref{thm3.7} and Theorem \ref{thm3.8}) and
\begin{equation}\label{eq5.25}
\left\|\Omega_{j} \left(t\right) \right\|_{W_1^1} \to 0 \; \left(j \to \infty\right),
\end{equation}
where
\begin{multline}\label{eq5.26}
\Omega_{j} \left(t\right):= \sum\limits_{n = 0}^\infty \left( \left| \left(\lambda_{n}-\lambda_{n,\left(j\right)} \right) \sin \lambda_{n} t  \right| \right.+\\
\left.+ \left| \left(\lambda_{n}-\lambda_{n,\left(j\right)} \right)^2 \cos \lambda_{n} t  \right| + \left| \left(k_{n}-k_{n,\left(j\right)} \right) \cos \lambda_{n} t  \right| \right),
\end{multline}
and $W_1^N$ the Sobolev space of functions $f\left(x\right),$ $x \in \left[0, \pi \right],$ such that $f^{\left(i\right)} \left(x \right),$ $i=0,1,\dots,N-1$ are absolutely continuous, equipped with the following norm:
\begin{equation}\label{eq5.27}
\| f\|_{W_1^N} = \sum_{i=0}^N \int_0^{\pi} \left| f^{(i)}(t) \right| dt.
\end{equation}
Denote
\begin{equation}\label{eq5.28}
H_j\left(t\right):=\sum\limits_{n = 0}^{\infty}\left(\dfrac{1}{a_{n,\left(j\right)}} \, \dfrac{\cos \lambda_{n,\left(j\right)} t}{\lambda_{n,\left(j\right)}^2}- \dfrac{1}{{a_n} \left(0, \pi, \beta \right)} \, \dfrac{\cos \lambda_{n}(0, \pi, \beta) t}{\lambda_{n}^{2}(0, \pi, \beta)}\right).
\end{equation}
By the same arguments as in Lemma \ref{lem4.2} one can prove that $H_{j}'\left(\cdot\right) \in AC\left(0, 2\pi \right).$ Let $P_j\left(x,t\right)$ be the solution of the following Gelfand-Levitan equation
\begin{equation*}
P_j \left(x, t \right)+F_j \left(x, t \right)+{\int_{0}^{x}} P_j \left(x, s \right)F_j \left(s, t \right)ds=0, \; 0 \leq t < x,
\end{equation*}
where $F_j \left(x, t \right)=\dfrac{1}{2} \left( H_j \left(\left|x-t \right| \right)-H_j \left(x+t \right) \right).$ Take
\begin{equation}\label{eq5.29}
\varphi_{\pi, j}\left(x, \mu\right):=\dfrac{\sin \lambda x}{\lambda} +{\displaystyle\int_{0}^{x}}P_j(x,t)\dfrac{\sin \lambda t} {\lambda}dt, \; q_j\left(x\right):=2\dfrac{d}{dx}P_j\left(x,x\right).
\end{equation}
Since $H_{j}'\left(\cdot\right) \in AC\left(0, 2\pi\right),$ then according to Case I, we receive that
\begin{equation}\label{eq5.30}
-\varphi_{\pi, j}''\left(x, \mu \right)+q_j\left(x\right)\varphi_{\pi, j} \left(x, \mu \right)=\mu \varphi_{\pi, j} \left(x, \mu \right),
\end{equation}
\begin{equation}\label{eq5.31}
\varphi_{\pi, j}\left(0, \mu \right)=0, \; \varphi_{\pi, j}'\left(0, \mu \right)=1.
\end{equation}
On one hand, according to \eqref{eq5.25} and Lemma \ref{lem7.2},
\begin{equation*}
\lim\limits_{j \to \infty}\| H_j(t)-H(t)\|_{W_1^1}=0,
\end{equation*}
which implies (by taking into account \cite[Lemma 1.5.1 on p. 32]{Freiling_Yurko:2001}) that
\begin{equation}\label{eq5.32}
\lim\limits_{j \to \infty} \mathop {\operatorname{max}} \limits_{0 \leq t \leq x \leq \pi}| P_j\left(x, t\right) - P\left(x, t\right)|=0
\end{equation}
and
\begin{equation}\label{eq5.33}
\lim\limits_{j \to \infty}\| q_j-q\|_{L_\mathbb{R}^1}=0.
\end{equation}
From \eqref{eq5.3}, \eqref{eq5.29} and \eqref{eq5.32}, it follows that
\begin{equation*}
\lim\limits_{j \to \infty} \mathop {\operatorname{max}} \limits_{0 \leq x \leq \pi} \mathop {\operatorname{max}} \limits_{|\mu| \leq r} | \varphi_{\pi, j}(x, \mu)-\varphi_{\pi}(x, \mu)|=0.
\end{equation*}
On the other hand, taking into account \cite[Theorem 3 on p. 62]{Harutyunyan_Hovsepyan:2005} (see, also, \cite[Theorem 5 on p. 18] {Poschel_Trubowitz:1987}), definition of the function $\tilde \varphi \left(x, \mu \right)$ and \eqref{eq5.33} we arrive at
\begin{equation*}
\lim\limits_{j \to \infty} \mathop {\operatorname{max}} \limits_{0 \leq x \leq \pi} \mathop {\operatorname{max}} \limits_{|\mu| \leq r} | \varphi_{\pi, j}(x, \mu)-\tilde \varphi(x, \mu)|=0.
\end{equation*}
Therefore $\tilde \varphi \left(x, \mu \right) \equiv \varphi_{\pi} \left(x, \mu \right).$ Lemma \ref{lem5.4} is proved.
\end{proof}
Using the same methods as in \cite{Gelfand_Levitan:1951,Gasymov_Levitan:1964,Freiling_Yurko:2001}, it can be proven that for arbitrary functions $f \in L^2 \left(0, \pi\right)$
\begin{equation}\label{eq5.34}
\int_{0}^{\pi} f^2\left(x \right)dx=\sum\limits_{n = 0}^\infty \dfrac{1}{a_n}\left(\int_{0}^{\pi} f\left(t \right)\varphi_{\pi}\left(t, \mu_n \right)dt\right)^2,
\end{equation}
or which is equivalent
\begin{equation}\label{eq5.35}
\int_{0}^{\pi} f\left(x \right)g\left(x \right)dx=\sum\limits_{n = 0}^\infty \dfrac{1}{a_n}\int_{0}^{\pi} f\left(t \right)\varphi_{\pi}\left(t, \mu_n \right)dt \int_{0}^{\pi} g\left(t \right)\varphi_{\pi}\left(t, \mu_n \right)dt,
\end{equation}
for every $f, g \in L^2 \left(0, \pi\right).$
\begin{lemma}\label{lem5.5}
The following relation holds
\begin{equation}\label{eq5.36}
\int_{0}^{\pi} \varphi_{\pi}\left(t,\mu_k\right)\varphi_{\pi}\left(t,\mu_n\right)dt = \begin{cases}
  0, \quad n \neq k, \\
  a_n, \quad n=k.
  \end{cases}
\end{equation}
\end{lemma}
\begin{proof}
Let $f' \in AC\left[0,\pi\right]$ and $f\left(0\right)=0$. Consider the following series
\begin{equation}\label{eq5.37}
f^{*}\left(x\right)=\sum\limits_{n = 0}^\infty c_n \varphi_{\pi}\left(x, \mu_n\right),
\end{equation}
where
\begin{equation}\label{eq5.38}
c_n :=\dfrac{1}{a_n}\int_{0}^{\pi} f\left(t\right)\varphi_{\pi} \left(t,\mu_n\right)dt.
\end{equation}
Using Lemma \ref{lem5.4} and integration by parts we obtain
\begin{multline}\label{eq5.39}
c_n=\dfrac{1}{a_n \mu_n}\int_{0}^{\pi} f\left(t\right)\left(-\varphi_{\pi}'' \left(t,\mu_n\right)+q\left(t\right)\varphi_{\pi} \left(t,\mu_n\right)\right)dt=\\
=\dfrac{1}{a_n \mu_n}\left(f'\left(\pi \right)\varphi_{\pi} \left(\pi,\mu_n\right)-f\left(\pi \right)\varphi_{\pi}' \left(\pi,\mu_n\right)\right)+\\
+\dfrac{1}{a_n \mu_n}\int_{0}^{\pi} \left(-f''\left(t\right)+q\left(t\right)f\left(t\right)\right)\varphi_{\pi} \left(t,\mu_n\right)dt.
\end{multline}
Since $\dfrac{1}{a_n \mu_n}=O\left(1\right)$ (see \eqref{eq1.11} and \eqref{eq1.12}) and
\begin{equation*}
\varphi_{\pi} \left(x, \mu_n \right)=\dfrac{\sin \left(n+\frac{1}{2}\right)x}{n+\frac{1}{2}}+O\left(\dfrac{1}{n^2}\right), \; \varphi_{\pi}' \left(x, \mu_n \right)=\cos \left(n+\frac{1}{2}\right)x+O\left(\dfrac{1}{n}\right)
\end{equation*}
uniformly on $\left[0, \pi \right],$ then $c_n=O\left(\dfrac{1}{n}\right).$ Hence the series \eqref{eq5.37} converges absolutely and uniformly on $\left[0, \pi \right].$ According to \eqref{eq5.35} and \eqref{eq5.38}, we have
\begin{multline*}
\int_{0}^{\pi} f\left(x \right)g\left(x \right)dx=\sum\limits_{n = 0}^\infty c_n \int_{0}^{\pi} g\left(t \right)\varphi_{\pi}\left(t, \mu_n \right)dt=\\
=\int_{0}^{\pi} g\left(t \right)\sum\limits_{n = 0}^\infty c_n \varphi_{\pi}\left(t, \mu_n \right)dt=\int_{0}^{\pi} g\left(t \right)f^{*}\left(t \right)dt.
\end{multline*}
As $g\left(x\right)$ is arbitrary, we get $f^{*}\left(x\right)=f\left(x\right),$ i.e.
\begin{equation*}
f\left(x\right)=\sum\limits_{n = 0}^\infty c_n \varphi_{\pi}\left(x, \mu_n\right).
\end{equation*}
Now, take $f\left(x\right)=\varphi_{\pi}\left(x, \mu_k\right)$ ($k \geq 0$ fixed). The system of functions $\left\{\varphi_{\pi}\left(x, \mu_n\right)\right\}_{n \geq 0}$ is minimal in $L^2\left(0, \pi \right)$ and therefore
\begin{equation*}
c_{nk} :=\dfrac{1}{a_n}\int_{0}^{\pi} \varphi_{\pi} \left(t,\mu_k\right)\varphi_{\pi} \left(t,\mu_n\right)dt=\delta_{nk},
\end{equation*}
where $\delta_{nk}$ is Kronecker symbol. Lemma \ref{lem5.5} is proved.
\end{proof}
\begin{lemma}\label{lem5.6}
For all $n, m \geq 0$
\begin{equation}\label{eq5.40}
\dfrac{\varphi_{\pi}'\left(\pi, \mu_n\right)}{\varphi_{\pi}\left(\pi, \mu_n\right)}=
\dfrac{\varphi_{\pi}'\left(\pi, \mu_m\right)}{\varphi_{\pi}\left(\pi, \mu_m\right)}=const
\end{equation}
\end{lemma}
\begin{proof}
It follows from \eqref{eq5.5} that
\begin{equation*}
\left.\left(\varphi_{\pi}\left(x, \mu_n\right)\varphi_{\pi}'\left(x, \mu_m\right)-\varphi_{\pi}'\left(x, \mu_n\right)\varphi_{\pi}\left(x, \mu_m\right)\right)\right|_{0}^{\pi}=\left(\mu_n-\mu_m\right)\int_{0}^{\pi} \varphi_{\pi}\left(x, \mu_n\right)\varphi_{\pi}\left(x, \mu_m\right)dx.
\end{equation*}
According to \eqref{eq5.36},
\begin{equation}\label{eq5.41}
\varphi_{\pi}\left(\pi, \mu_n\right)\varphi_{\pi}'\left(\pi, \mu_m\right)-\varphi_{\pi}'\left(\pi, \mu_n\right)\varphi_{\pi}\left(\pi, \mu_m\right)=0.
\end{equation}
Let us show that $\varphi_{\pi}\left(\pi, \mu_n\right) \neq 0$ for all $n \geq 0.$ Otherwise from $\varphi_{\pi}\left(\pi, \mu_m\right)=0$ for a certain $m$ and since $\varphi_{\pi}'\left(\pi, \mu_m\right) \neq 0,$ from \eqref{eq5.41} we get that $\varphi_{\pi}\left(\pi, \mu_n\right)=0$ for all $n,$ which is impossible since
\begin{equation*}
\left(n+\dfrac{1}{2}\right)\varphi_{\pi}\left(\pi, \mu_n\right)=(-1)^n+O\left(\dfrac{1}{n}\right).
\end{equation*}
Then, dividing \eqref{eq5.41} by $\varphi_{\pi}\left(\pi, \mu_n\right)\varphi_{\pi}\left(\pi, \mu_m\right),$ we arrive at \eqref{eq5.40}. Lemma \ref{lem5.6} is proved.
\end{proof}
Denote
\begin{equation}\label{eq5.42}
\cot\tilde \beta := -\dfrac{\varphi_{\pi}'\left(\pi, \mu_n\right)}{\varphi_{\pi}\left(\pi, \mu_n\right)}.
\end{equation}
Thus,
\begin{equation*}
\varphi_{\pi}\left(\pi, \mu_n\right)\cos \tilde \beta+\varphi_{\pi}'\left(\pi, \mu_n\right)\sin \tilde \beta=0, n \geq 0.
\end{equation*}
Together with Lemma \ref{lem5.4} and Lemma \ref{lem5.5} this gives that $\left\{\mu_n\right\}_{n \geq 0}$ and  $\left\{a_n\right\}_{n \geq 0}$ are the eigenvalues and the norming constants for the constructed $L\left(q, \pi, \tilde \beta\right)$ problem, respectively.

\begin{remark}\label{rem5.7}
Note that, on one hand (see \eqref{eq1.11}, \eqref{eq3.1} and \eqref{eq3.8})
\begin{equation*}
\cot \tilde \beta = \cot \beta +\dfrac{1}{2} \left(\pi c - \displaystyle\int_{0}^{\pi} q \left( t \right) dt \right),
\end{equation*}
where $\tilde \beta$ and $q$ are defined in \eqref{eq5.42} and \eqref{eq5.4}, respectively and $\beta$ and $c$ determine the sequences $\left\{\mu_n \right\}_{n=0}^{\infty}$ and $\left\{a_n \right\}_{n=0}^{\infty}$ (see \eqref{eq1.11}, \eqref{eq1.12}).

On the other hand, it is easy to see that $\tilde \beta$ may not coincide with $\beta$ (see section \ref{sec6}). In \cite{Ashrafyan_Harutyunyan:2015} the authors proposed additional conditions (necessary and sufficient) that ensure $\alpha=\tilde \alpha$ and $\beta=\tilde \beta,$ when $q \in L_{\mathbb{R}}^2\left[ 0, \pi \right],$ $\alpha, \beta \in \left(0, \pi \right)$ problem is considered. We believe that analogous condition can be found in our case as well, but this is outside the scope of the present paper.
\end{remark}
The Theorem \ref{thm1.1} is completely proved.

\section{Implementation of the algorithm}
\label{sec6}

In this section we are going to present an example of realization of the constructive solution of an inverse Sturm-Liouville problem. To this aim, we will take two sequences $\left\{\mu_n\right\}_{n=0}^{\infty}$ and $\left\{a_n\right\}_{n=0}^{\infty}$ satisfying the conditions of the Theorem \ref{thm1.1} and applying the procedure describing in the Section \ref{sec5} to reconstruct the potential $q$ and the parameter $\tilde \beta.$ Let
\begin{equation}\label{eq6.1}
\sqrt {\mu_n} \equiv \lambda_n = n+\dfrac{1}{2}, \; n=0,1,2,\dots,
\end{equation}
\begin{equation}\label{eq6.2}
a_0=\pi, \; a_n = \dfrac{\pi}{2 \left(n+\frac{1}{2}\right)^2}, \; n=1,2,\dots,
\end{equation}
then $\beta=\dfrac{\pi}{2}$ (see \eqref{eq1.11}, \eqref{eq1.12} and \eqref{eq3.8}) and according to \eqref{eq4.15}
\begin{multline}\label{eq6.3}
F\left(x, t \right)=\sum\limits_{n=0}^{\infty}\left(\dfrac{1}{a_n} \, \dfrac{\sin \lambda_{n} x}{\lambda_{n}} \, \dfrac{\sin \lambda_{n} t}{\lambda_{n}}- \dfrac{2}{\pi} \, {\sin \left(n+\dfrac{1}{2} \right)x} \, {\sin \left(n+\dfrac{1}{2} \right)t} \right)=\\
=\dfrac{2}{\pi} \, \sin \dfrac{x}{2} \, \sin \dfrac{t}{2}.
\end{multline}
We will seek the solution $P\left(x, t\right)$ of the corresponding Gelfand-Levitan equation (see \eqref{eq4.16}) in the form $P\left(x, t\right)=a\left(x\right)\sin \dfrac{t}{2}.$ After some calculations we find that
\begin{equation}\label{eq6.4}
P\left(x, t\right)=\dfrac{4}{2 \sin x -2x-2\pi} \, \sin \dfrac{x}{2} \, \sin \dfrac{t}{2}
\end{equation}
and
\begin{equation}\label{eq6.5}
q\left(x \right):=2\dfrac{d}{dx}P\left(x, x \right)=\dfrac{2 \sin x}{\sin x -x-\pi}-\dfrac{4 \left(\cos x -1\right)}{\left(\sin x -x-\pi\right)^2} \, \sin^2 \dfrac{x}{2},
\end{equation}
\begin{equation}\label{eq6.6}
\cot\tilde \beta = -\dfrac{\varphi_{\pi}'\left(\pi, \mu_n\right)}{\varphi_{\pi}\left(\pi, \mu_n\right)}=-P \left(\pi, \pi \right)=-\dfrac{1}{\pi},
\end{equation}
\begin{equation}\label{eq6.7}
\tilde \beta = arccot \dfrac{1}{\pi}.
\end{equation}

\section{Appendix}
\label{app}

Here, we prove two lemmas, which play an important role in our analysis.
\begin{lemma}\label{lem7.1}
Let us denote
\begin{equation}\label{eq7.1}
T_\beta \left(x \right):= \sum_{n=2}^{\infty} \dfrac{\sin \left(n+\delta_n\left(\pi, \beta\right)\right)x}{n+\delta_n\left(\pi, \beta\right)}.
\end{equation}
Then $T_\beta \left(x \right)$ and $T'_\beta \left(x \right)$ are absolutely continuous functions on arbitrary segment $\left[a, b\right] \subset \left(0, 2\pi\right).$
\end{lemma}
\begin{proof}
Denote $t_n=\delta_n\left(\pi, \beta\right)-\dfrac{1}{2}$ and write the general term of the series in \eqref{eq7.1} in the following form
\begin{multline}\label{eq7.2}
\dfrac{\sin \left(n+\delta_n\left(\pi, \beta\right)\right)x}{n+\delta_n\left(\pi, \beta\right)}=\dfrac{\sin \left(n+\frac{1}{2}+t_n \right)x}{n+\frac{1}{2}}-\dfrac{t_n\sin \left(n+\frac{1}{2}+t_n\right)x}{\left(n+\frac{1}{2}\right)\left(n+\frac{1}{2}+t_n\right)}=\\
=\dfrac{\sin \left(n+\frac{1}{2}\right)x}{n+\frac{1}{2}} \, \cos t_{n}x + \dfrac{\cos \left(n+\frac{1}{2}\right)x}{n+\frac{1}{2}} \, \sin t_{n}x - \dfrac{t_n\sin \left(n+\frac{1}{2}\right)x}{\left(n+\frac{1}{2}\right)\left(n+\frac{1}{2}+t_n\right)}
\end{multline}
It follows from \eqref{eq3.8}, that
\begin{equation*}
t_n=\dfrac{\cot \beta}{\pi \left(n+\frac{1}{2}\right)}+u_n, \quad \cos t_{n}x=1+v_n \left(x \right), \quad \sin t_{n}x=t_{n}x+w_n \left(x \right),
\end{equation*}
where $u_n=O\left(\dfrac{1}{n^2}\right)$ and $v_n \left(x \right)=O\left(\dfrac{1}{n^2}\right),$ $w_n \left(x \right)=O\left(\dfrac{1}{n^3}\right)$ uniformly on $\left[0, 2\pi\right]$ and therefore \eqref{eq7.2} can be written in the following way:
\begin{multline*}
\dfrac{\sin \left(n+\delta_n\left(\pi, \beta\right)\right)x}{n+\delta_n\left(\pi, \beta\right)}=\dfrac{\sin \left(n+\frac{1}{2}\right)x}{n+\frac{1}{2}}+v_n \left(x \right) \dfrac{\sin \left(n+\frac{1}{2}\right)x}{n+\frac{1}{2}}+\dfrac{x \cot\beta}{\pi} \, \dfrac{\cos \left(n+\frac{1}{2}\right)x}{\left(n+\frac{1}{2}\right)^2}+\\
+x \, u_{n} \, \dfrac{\cos \left(n+\frac{1}{2}\right)x}{n+\frac{1}{2}}+ w_n \left(x \right) \dfrac{\cos \left(n+\frac{1}{2}\right)x}{n+\frac{1}{2}}-\dfrac{t_n\sin \left(n+\frac{1}{2}\right)x}{\left(n+\frac{1}{2}\right)\left(n+\frac{1}{2}+t_n\right)}.
\end{multline*}
It is easy to see that
\begin{equation*}
T_{\beta} \left( x \right)=T_{1} \left( x \right)+T_{2\beta} \left( x \right)+T_{3\beta} \left( x \right),
\end{equation*}
where
\begin{equation*}
T_{1} \left( x \right)=\sum_{n=2}^{\infty}\dfrac{\sin \left(n+\frac{1}{2}\right)x}{n+\frac{1}{2}},
\end{equation*}
\begin{equation*}
T_{2\beta} \left( x \right)=\dfrac{x \cot\beta}{\pi}\sum_{n=2}^{\infty} \dfrac{\cos \left(n+\frac{1}{2}\right)x}{\left(n+\frac{1}{2}\right)^2},
\end{equation*}
\begin{multline*}
T_{3\beta} \left( x \right)=\sum_{n=2}^{\infty} \left(v_n \left(x \right) \dfrac{\sin \left(n+\frac{1}{2}\right)x}{n+\frac{1}{2}}+x \, u_{n} \, \dfrac{\cos \left(n+\frac{1}{2}\right)x}{n+\frac{1}{2}}+ \right. \\
\left. +w_n \left(x \right) \dfrac{\cos \left(n+\frac{1}{2}\right)x}{n+\frac{1}{2}}-\dfrac{t_n\sin \left(n+\frac{1}{2}\right)x}{\left(n+\frac{1}{2}\right)\left(n+\frac{1}{2}+t_n\right)}\right).
\end{multline*}
From \cite[formulae 37 and 38 on p. 578]{Bronshtein_Semendyayev:1998}, we get that
\begin{equation*}
T_{1} \left( x \right)=\dfrac{\pi}{2}-2\sin\dfrac{x}{2}-\dfrac{2}{3}\sin\dfrac{3x}{2}, \quad 0 < x < 2\pi
\end{equation*}
and
\begin{equation*}
T_{2\beta} \left( x \right)=\dfrac{x \cot\beta}{\pi} \left(\dfrac{\pi^2-2\pi x}{2}-4\cos \dfrac{x}{2}-\dfrac{4}{9}\cos\dfrac{3x}{2} \right), \quad 0 \leq x \leq 2\pi
\end{equation*}
are infinitely differentiable functions on corresponding domains. The estimates for $u_n,$ $v_n \left(x \right)$ and $w_n \left(x \right)$ ensure the absolutely continuity of the functions $T_{3\beta}$ and $T'_{3\beta}.$ This completes the proof.
\end{proof}

\begin{lemma}\label{lem7.2}
Let us given the functions $\Omega_{j} \left(t\right),$ $H_j\left(t \right)$ and $H\left(t \right)$ defined in \eqref{eq5.26}, \eqref{eq5.28} and \eqref{eq4.3}, respectively. Then, there are $C_1$ and $C_2$ positive numbers such that
\begin{equation}\label{eq7.3}
\mathop {\operatorname{max}} \limits_{0 \leq t \leq 2\pi}\left| H_j\left(t \right) - H\left(t \right) \right| \leq C_1
\mathop {\operatorname{max}} \limits_{0 \leq t \leq 2\pi}\left| \Omega_{j} \left(t\right) \right|,
\end{equation}
\begin{equation}\label{eq7.4}
\left\|H_j \left(t\right) - H \left(t\right) \right\|_{W_1^1} \leq C_2 \left\| \Omega_{j} \left(t\right) \right\|_{W_1^1}.
\end{equation}
\end{lemma}

\begin{proof}
Using \eqref{eq5.28} and \eqref{eq4.3} we write down $H_j\left(t \right) - H\left(t \right)$ and $H_j'\left(t \right) - H'\left(t \right)$ in the following forms:
\begin{multline}\label{eq7.5}
H_j\left(t \right) - H\left(t \right) = \\
\sum\limits_{n = 0}^\infty \left( \dfrac{1}{a_{n,\left(j\right)}} \, \dfrac{\cos \lambda_{n,\left(j\right)} t}{\lambda_{n,\left(j\right)}^2}- \dfrac{1}{a_n} \, \dfrac{\cos \lambda_{n} t}{\lambda_{n}^2} \right)= \sum\limits_{n = 0}^\infty \left( \dfrac{\cos \lambda_{n,\left(j\right)} t}{k_{n,\left(j\right)}}- \dfrac{{\cos \lambda_{n} t}}{k_n} \right)= \\
=\sum\limits_{n = 0}^\infty \left( \dfrac{1}{k_{n,\left(j\right)}} \left(\cos \lambda_{n,\left(j\right)} t-\cos \lambda_{n} t\right)+\left(\dfrac{1}{k_{n,\left(j\right)}}-\dfrac{1}{k_{n}}\right) {\cos \lambda_{n} t}\right)=\\
=\sum\limits_{n = 0}^\infty \left( \dfrac{2}{k_{n,\left(j\right)}} \sin \dfrac{\left(\lambda_{n}-\lambda_{n,\left(j\right)}\right)t}{2} \, \sin \dfrac{\left(\lambda_{n}+\lambda_{n,\left(j\right)}\right)t}{2} +\dfrac{k_{n}-k_{n,\left(j\right)}}{k_{n} \, k_{n,\left(j\right)}} {\cos \lambda_{n} t}\right)=\\
= \sum\limits_{n = 0}^\infty \left( \dfrac{2}{k_{n,\left(j\right)}} \sin \dfrac{\left(\lambda_{n}-\lambda_{n,\left(j\right)}\right)t}{2} \, \sin \left(\lambda_{n} t-  \dfrac{\left(\lambda_{n}-\lambda_{n,\left(j\right)}\right)t}{2}\right) +\right.\\
\left.+\dfrac{k_{n}-k_{n,\left(j\right)}}{k_{n} \, k_{n,\left(j\right)}} {\cos \lambda_{n} t}\right)=
\sum\limits_{n = 0}^\infty \left( \dfrac{1}{k_{n,\left(j\right)}} \sin \left(\lambda_{n}-\lambda_{n,\left(j\right)}\right)t \, \sin \lambda_{n} t- \right. \\
\left.-\dfrac{2}{k_{n,\left(j\right)}} \sin^2 \dfrac{\left(\lambda_{n}-\lambda_{n,\left(j\right)}\right)t}{2} \, \cos \lambda_{n} t +\dfrac{k_{n}-k_{n,\left(j\right)}}{k_{n} \, k_{n,\left(j\right)}} {\cos \lambda_{n} t}\right),
\end{multline}
\begin{multline}\label{eq7.6}
H_j'\left(t \right) - H'\left(t \right) =\sum\limits_{n = 0}^\infty \left( \dfrac{\lambda_{n}-\lambda_{n,\left(j\right)}}{k_{n,\left(j\right)}} \cos \left(\lambda_{n}-\lambda_{n,\left(j\right)}\right)t \, \sin \lambda_{n} t+\right. \\
\left.+\dfrac{\lambda_{n}}{k_{n,\left(j\right)}} \sin \left(\lambda_{n}-\lambda_{n,\left(j\right)}\right)t \, \cos \lambda_{n} t\right)-\sum\limits_{n = 0}^\infty \left( \dfrac{\lambda_{n}-\lambda_{n,\left(j\right)}}{k_{n,\left(j\right)}} \sin \left(\lambda_{n}-\lambda_{n,\left(j\right)}\right)t \, \cos \lambda_{n} t -\right. \\
\left.-\dfrac{2 \lambda_{n}}{k_{n,\left(j\right)}} \sin^2 \dfrac{\left(\lambda_{n}-\lambda_{n,\left(j\right)}\right)t}{2} \, \sin \lambda_{n} t\right)- \sum\limits_{n = 0}^\infty \dfrac{k_n-k_{n,\left(j\right)}}{k_n k_{n,\left(j\right)}} \lambda_{n} \sin \lambda_{n} t.
\end{multline}
From \eqref{eq7.5} and \eqref{eq7.6} it is easy to get the following estimates:
\begin{equation}\label{eq7.7}
\left|H_j\left(t \right) - H\left(t \right)\right| \leq C_1 \Omega_{j} \left(t\right),
\end{equation}
\begin{equation}\label{eq7.8}
\left|H_j'\left(t \right) - H'\left(t \right)\right| \leq C_1 \left( \Omega_{j} \left(t\right) + \Omega_{j}' \left(t\right)\right).
\end{equation}
These prove the lemma.
\end{proof}

\subsection*{Acknowledgment}
The author was supported by State Committee of Science MES RA, in frame of the research project No. 16YR--1A017. The author is grateful to his supervisor, Professor Tigran Harutyunyan for suggesting the problem and attention to the work, anonymous referee whose valuable suggestions and comments improve the first version of this manuscript and Dr. Avetik Arakelyan for carefully reading the article and valuable remarks.


\begin{thebibliography}{32}

\bibitem{Ashrafyan_Harutyunyan:2015}
Ashrafyan, Yu.A. and Harutyunyan, T.N.
``Inverse Sturm-Liouville problems with fixed boundary conditions.''
\textit{Electron. J. Diff. Equ.}, 2015, no. 27, (2015): 1--8.

\bibitem{Atkinson:1964}
Atkinson, F.V.
\textit{Discrete and continuous boundary problems},
Academic Press, New York-London, 1964.

\bibitem{Bari:1961}
Bari, N.K.
\textit{Trigonometric Series},
Fizmatgiz, Moscow, (in Russian), 1961.

\bibitem{Bronshtein_Semendyayev:1998}
Bronshtein, I.N. and Semendyayev K.A.
\textit{Handbook of Mathematics},
Springer-Verlag Berlin Heidelberg, New York, 1998.

\bibitem{Chudov:1949}
Chudov, L.A.
``The inverse Sturm-Liouville problem.''
\textit{Mat. Sbornik}, 25(67), no. 3, (in Russian), (1949): 451--456.

\bibitem{Coddington_Levinson:1955}
Coddington, E. and Levinson N.
\textit{Theory of Ordinary Differential Equations},
McGraw Hill Book Company, New York, 1955.

\bibitem{Dahlberg_Trubozitz:1984}
Dahlberg, B.E.J. and Trubowitz, E.
``The inverse Sturm-Liouville problem. III.''
\textit{Comm. Pure Appl. Math.}, 37, no. 2, (1984): 255--267.

\bibitem{Freiling_Yurko:2001}
Freiling, G. and Yurko, V.A.
\textit{Inverse Sturm-Liouville Problem and Their Applications},
NOVA Science Publishers, New York, 2001.

\bibitem{Gasymov_Levitan:1964}
Gasymov, M.G. and Levitan, B.M.
``Determination of a differential equation by two of its spectra.''
\textit{Uspekhi Mat. Nauk}, 19, no. 2, (in Russian), (1964): 3--63.

\bibitem{Gelfand_Levitan:1951}
Gelfand, I.M. and Levitan, B.M.
``On the determination of a differential equation from its spectral function.''
\textit{Izv. Akad. Nauk SSSR, ser. Math.}, 15, no. 4, (in Russian), (1951): 309--360.

\bibitem{Harutyunyan:2008}
Harutyunyan, T.N.
``The Dependence of the Eigenvalues of the Sturm-Liouville Problem on Boundary Conditions.''
\textit{Matematicki Vesnik}, 60, no. 4, (2008): 285--294.

\bibitem{Harutyunyan:2010}
Harutyunyan, T.N.
\textit{Eigenvalue Functions of Family of Sturm-Liouville and Dirac Operators},
Doctoral Thesis, Yerevan, (in Russian), 2010.

\bibitem{Harutyunyan:2016}
Harutyunyan, T.N.
``Asymptotics of the eigenvalues of Sturm-Liouville problem.''
\textit{Journal of Contemporary Mathematical Analysis}, 51, no. 4, (2016): 174--183.

\bibitem{Harutyunyan_Hovsepyan:2005}
Harutyunyan, T.N. and Hovsepyan, M.S.
``On the solutions of the Sturm-Liouville equation.''
\textit{Mathematics in Higher School}, 1, no. 3, (in Russian), (2005): 59--74.

\bibitem{Harutyunyan_Pahlevanyan:2016}
Harutyunyan, T.N. and Pahlevanyan, A.A.
``On the norming constants of the Sturm-Liouville problem.''
\textit{Bulletin of Kazan State Power Engineering University}, no. 3(31), (2016): 7--26.

\bibitem{Harutyunyan_Pahlevanyan_Srapionyan:2013}
Harutyunyan, T.N., Pahlevanyan, A.A. and Srapionyan A.V.
``Riesz bases generated by the spectra of Sturm-Liouville problem.''
\textit{Electron. J. Diff. Equ.}, 2013, no. 71, (2013): 1--8.

\bibitem{Isaacson_Mckean_Trubowitz:1984}
Isaacson, E.L., McKean, H.P. and Trubowitz, E.
``The inverse Sturm-Liouville problem. II.''
\textit{Comm. Pure Appl. Math.}, 37, no. 1, (1984): 1--11.

\bibitem{Isaacson_Trubowitz:1983}
Isaacson, E.L. and Trubowitz, E.
``The inverse Sturm-Liouville problem. I.''
\textit{Comm. Pure Appl. Math.}, 36, no. 6, (1983): 767--783.

\bibitem{Iserles_Norsett:2008}
Iserles, A. and N{\o}rsett, S.P.
``From high oscillation to rapid approximation I: Modified Fourier expansions.''
\textit{IMA Journal of Numerical analysis}, 28, no. 4, (2008): 862--887.

\bibitem{Kadec:1964}
Kadec, M.I.
``The exact value of the Paley-Wiener constant.''
\textit{Sov. Math Doklady}, 5, no. 2, (1964): 559--561.

\bibitem{Korotyaev_Chelkak:2009}
Korotyaev, E.L. and Chelkak D.S.
``The inverse Sturm–Liouville problem with mixed boundary conditions.''
\textit{Algebra i Analiz}, 21, no. 5, (in Russian), (2009): 114-–137.

\bibitem{Levinson:1940}
Levinson, N.
\textit{Gap and density theorems},
American Mathematical Society, New York, 1940.

\bibitem{Levitan:1984}
Levitan, B.M.
\textit{Inverse Sturm-Liouville problems},
Nauka, Moscow, (in Russian), 1984.

\bibitem{Levitan_Sargsyan:1970}
Levitan, B.M. and Sargsyan I.S.
\textit{Introduction to Spectral Theory},
Nauka, Moscow, (in Russian), 1970.

\bibitem{Levitan_Sargsyan:1988}
Levitan, B.M. and Sargsyan, I.S.
\textit{Sturm-Liouville and Dirac operators},
Nauka, Moscow, (in Russian), 1988.

\bibitem{Marchenko:1952}
Marchenko, V.A.
``Some questions of the theory of one-dimensional linear differential operators of the second order.''
\textit{Trudy Moskov. Mat. Obsh.}, 1, (in Russian), (1952): 327--420.

\bibitem{Marchenko:1977}
Marchenko, V.A.
\textit{The Sturm-Liouville Operators and their Applications},
Naukova Dumka, Kiev, (in Russian), 1977.

\bibitem{Mihlin:1949}
Mihlin, S.G.
\textit{Integral equations and their applications to some problems of mechanics, mathematical physics and engineering},
GITTL, Moscow-Leningrad, (in Russian), 1949.

\bibitem{Naimark:1969}
Naimark, M.A.
\textit{Linear Differential Equations},
Nauka, Moscow, (in Russian), 1969.

\bibitem{Poschel_Trubowitz:1987}
P{\"o}schel, J. and Trubowitz, E.
\textit{Inverse spectral theory},
Academic Press, Inc., Boston, MA, 1987.

\bibitem{Povzner:1948}
Povzner, A.
``On differential equations of Sturm-Liouville type on a half-axis.''
\textit{Mat. Sbornik}, 23(65), no. 1, (in Russian), (1948): 3--52.

\bibitem{Shabat:1985}
Shabat, B.V.
\textit{Introduction to complex analysis},
Nauka, Moscow, (in Russian), 1985.

\bibitem{Young:1980}
Young, R.M.
\textit{An introduction to nonharmonic Fourier series},
Academic Press, New York, 1980.

\bibitem{Zhikov:1967}
Zhikov, V.V.
``On inverse Sturm-Liouville problems on a finite segment.''
\textit{Izv. Akad. Nauk SSSR, ser. Math.}, 31, no. 5, (in Russian), (1967): 965--976.

\end{thebibliography}
\end{document}